%% file: 905.tex
\documentclass[a4paper]{amsart}
\usepackage{amsfonts,amssymb,amsmath,amsthm}
\usepackage{accents}
\usepackage{type1cm} 

   \usepackage[dvips]{hyperref}
   \usepackage[dvips]{graphicx}

\usepackage{color} 

\numberwithin{equation}{section}

\newtheorem{MainThm}{Theorem}
\newtheorem{Thm}{Theorem}[section]
\newtheorem{Lem}[Thm]{Lemma}
\newtheorem{Cor}[Thm]{Corollary}

\theoremstyle{definition}
\newtheorem{notation}[Thm]{Notation}
\newtheorem{Def}[Thm]{Definition}

\theoremstyle{remark}
\newtheorem{remark}[Thm]{Remark}
\newtheorem{remarks}[Thm]{Remarks}

\newtheorem{Fact}[Thm]{Fact}

\newenvironment{eqpar}{\begin{equation}\begin{minipage}{0.8\columnwidth}}%
{\end{minipage}\end{equation}}

\renewcommand{\restriction}{{\upharpoonright}}

\DeclareMathOperator{\ro}{ro}
\DeclareMathOperator{\dom}{dom}               
\DeclareMathOperator{\height}{ht}

\newcommand{\al}[1]{\ensuremath{{\aleph_{#1}}} }          
\newcommand{\om}[1]{\ensuremath{{\omega_{#1}}} }          

\newcommand{\esm}{\ensuremath{\prec}}  

\newcommand{\incomp}{\ensuremath{\perp}}

\newcommand{\n}[1]{\underaccent{\tilde}{#1}}

\newcommand{\rdelta}{r^{\Delta}}
\newcommand{\rmulti}{r^{\text{mult}}}
\begin{document}

\date{2007-03-11}

\title{A Sacks Real out of Nowhere}
\author[Jakob Kellner]{Jakob Kellner$^\dag$}
\address{Jakob Kellner\\Kurt G\"odel Research Center for Mathematical Logic
at the University of Vienna\\ W\"ahringer Stra\ss e 25\\ 1090 Wien, Austria}
\email{kellner@fsmat.at}
\urladdr{http://www.logic.univie.ac.at/$\sim$kellner}
\thanks{$^\dag$ supported by a European Union Marie Curie EIF fellowship, contract MEIF-CT-2006-024483.}
\author[Saharon Shelah]{Saharon Shelah$^\ddag$}
\address{Saharon Shelah\\Einstein Institute of Mathematics\\
Edmond J. Safra Campus, Givat Ram\\
The Hebrew University of Jerusalem\\
Jerusalem, 91904, Israel\\
and
Department of Mathematics\\
Rutgers University\\
New Brunswick, NJ 08854, USA}
\email{shelah@math.huji.ac.il}
\urladdr{http://www.math.rutgers.edu/$\sim$shelah}
\thanks{
$^\ddag$
supported by the United States-Israel
  Binational Science Foundation (Grant no. 2002323),
publication 905.}

\begin{abstract}
  There is a proper countable support iteration
  of length $\omega$ adding no new reals at finite stages
  and adding a Sacks real in the limit.
\end{abstract}

\maketitle

\section{Introduction}

Preservation theorems are a central tool in forcing theory:
\begin{quote}
  Let $(P_\alpha,\n Q_\alpha)_{\alpha<\epsilon}$ be a forcing iteration. Assume
  that $Q_\alpha$ is (forced to be) 
  nice for all $\alpha<\epsilon$. Then $P_\epsilon$ is nice.\footnote{Often
    preservation theorems also have the (weaker) form ``If 
    $\epsilon$ is a limit, and all $P_\alpha$ are nice for $\alpha<\epsilon$,
    then $P_\epsilon$ is nice.''}
\end{quote}
A niceness (or preservation) property usually implies that the forcing does
not change the universe too much.
Among the most important preservation theorems are:
\begin{center}
  The finite support iteration of ccc forcings is ccc. \cite{MR0294139}
\end{center}
and
\begin{center}
  The countable support iteration of proper forcings is proper. \cite{MR675955}
\end{center}
In this paper we investigate proper countable support iterations, so the
limits are always proper. Many additional preservation properties are
preserved as
well, for example $\omega^\omega$-bounding (i.e., not adding an unbounded real).
This is a special instance of a general preservation theorem by the second
author (``Case A'' of \cite[XVIII \S3]{MR1623206}) 
which is also known as ``first preservation theorem'' \cite[Section
6.1.B]{MR1350295}
or ``tools-preservation'' \cite[Section 5]{MR1234283}, see also
\cite[Theorem 2.4]{MR2214624}.
Many additional preservation theorems for proper countable support iterations
can be found in \cite{MR1623206}, or,
from the point of view of large cardinals, in \cite{MR2391923}.

We investigate iterations where all iterands are NNR, which means that they do
not add new reals. So the iterands (and therefore the limit as well) satisfy
all instances
of tools-preservation. 
However, it turns out that 
the limit can add a new real $r$.
The first example was given by Jensen \cite{MR0384542}, and the 
phenomenon was further investigated in \cite[V]{MR1623206}.
So what do we know about the real $r$?
We know that it has to be bounded by an old real (i.e., a real in the ground
model), corresponding to the iterable preservation property
``$\omega^\omega$-bounding''.  $r$ will even satisfy the stronger Sacks
property.  In particular $r$ cannot be, e.g., a Cohen, random, Laver or Mathias real.
In the previously known examples, 
the proof that a new real $r$ is added is rather
indirect and does not give much ``positive'' information about $r$.
So it is natural to ask which kind of reals
{\em can} appear in proper NNR  limits.
Todd Eisworth asked this question for the simplest and best understood real
that satisfies the Sacks property, the Sacks real. In this paper, we show that
Sacks reals indeed can appear in this way:
\begin{MainThm}\label{thm:main}
  There is an iteration $(P_n,\n Q_n)_{n<\omega}$
  such that each $Q_n$ is forced to be proper and NNR
  and such that the countable support limit
  $P_\omega$ adds a Sacks real.
  Moreover, $P_\omega$ is equivalent to $S\ast P'$, where 
  $S$ is Sacks forcing and $P'$ is NNR.\footnote{We 
    do not claim that $P'$ is proper.}
\end{MainThm}

The Theorem can be interpreted in two ways:

On the one hand, it indicates limitations of possible preservation theorems:
``Not adding a Sacks real'' is obviously not iterable (even with rather strong
additional assumptions).

On the other hand, it shows that Sacks forcing is exceptionally ``harmless'':
It satisfies every usual iterable preservation property.\footnote{More exactly:
  Sacks forcing satisfies every property $X$ such that:
  Every proper NNR forcing satisfies $X$, 
  $X$ is preserved under proper countable support iterations,
  and if $P$ does {\em not} satisfy $X$, then $P\ast \n Q$ does not satisfy
  $X$ either for any NNR $Q$.}
So the Sacks model (the model constructed by starting with CH and iterating
$\omega_2$ Sacks forcings in a countable support iteration) has all the
corresponding properties as well.\footnote{Of course this is already known for
  many of the popular properties, cf.~ \cite{MR2176267} or
  \cite{MR2391923}, which shows that in some respect Sacks forcing is the ``most
  tame'' forcing possible. This corresponds to the fact that all of the usual
  cardinal characteristics (apart from the continuum) are $\al1$ in the Sacks
  model.}

In a continuation of this work we will say more about the kind of reals that
can be added in limits of NNR iterations (e.g. generics for other finite
splitting lim sup tree forcings). It turns out that  many of these reals can
appear at limit stages, but some of them not at stage $\omega$, but only at
later stages, e.g.\ $\omega^2$.

We thank a referee for suggesting several improvements in the presentation.

\section{Sacks conditions as squares of terms}\label{sec:sacks1}

In this section, we introduce the forcing notion $Q_*$, which 
is forcing equivalent
to Sacks forcing. We will later work with $Q_*$ in the proof of 
Theorem~\ref{thm:main}.

A Sacks condition (or Sacks tree) is a perfect tree $T\subseteq 2^{<\omega}$.
Given $T$, we call a node  $t$ a splitting node if $t$ has two immediate
successors in $T$.

Let $F_n$ be the set of the $n$-th splitting nodes, cf.~Figure~\ref{fig:tree}.
So 
$t\in F_n$ means that $t$ is a splitting node and that there are $n$ splitting
nodes below $t$.  Since $T$ is perfect, $F_n$ is a front, which means that
every branch through $T$ meets $F_n$ exactly once.
Being a front is stronger than just being
a maximal antichain, and due to K\"onig's Lemma every front is finite.

A branch $b$ through $T$ is an element of $2^\omega$ and therefore a sequence
$(b_0,b_1,\dots)$ for some $b_n\in \{0,1\}$. Intuitively speaking, we can
describe ``the arbitrary branch'' $b$ of $T$ by interpreting each $b_n$ to be a
term $t_n(x_0,\dots,x_n)$, where the value of $t_n$ (0 or 1) depends on 
$x_l$ for $l\leq n$, and $x_l$ is a variable with values in $\{0,1\}$ that
tells us whether we choose the left (0) or right (1) path at the front $F_l$. 

A more formal description of terms can be found in Definition~\ref{def:terms},
but a simple example is much more instructive:
\begin{figure}[tb]
      \newcommand{\mylabeltree}{\relax}
    \begin{center}
      \scalebox{0.6}{\input{tree.pstex_t}}
    \end{center}
  \caption{\label{fig:tree}
  $F_n$ is the front of $n$-th splitting nodes. 
  }
\end{figure}
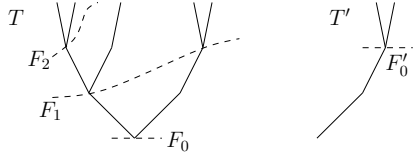
In the tree $T$ of Figure~\ref{fig:tree}, 
the sequence of terms begins as follows: 
\[
  t_0=x_0,
  \quad
  t_1=\begin{cases}
         x_1&\text{if }x_0=0,\\
         1  &\text{otherwise,}
       \end{cases}
  \quad
  t_2=\begin{cases}
         x_2&\text{if }x_0=0\text{ and }x_1=0,\\
         1  &\text{if }x_0=0\text{ and }x_1=1,\\
         x_1&\text{otherwise.}
       \end{cases}
\]

We will use the following notation:
\begin{eqpar}\label{eq:canonical}
  Given a Sacks tree $T$, the sequence $\bar t$ of terms defined as above is 
  called the canonical term sequence for $T$.
\end{eqpar}

Let $a$ be an assignment, that is a map that assigns each variable a value in
$\{0,1\}$. Then $a$ can be extended to evaluate terms $t$ to $t\circ
a\in\{0,1\}$, so we can evaluate the 
term sequence $\bar t=(t_0,t_1,\dots)$  to
\[
  \bar t\circ a:=(t_0\circ a,t_1\circ a,\dots)\in 2^\omega.
\]
If $\bar t$ is the  canonical term sequence for $T$ and $a$ an assignment, then
$\bar t\circ a$ is a branch through $T$.
Moreover, every branch can be obtained this way:
\begin{equation}\label{eq:blub2}
  T=\{\bar t \circ a\restriction n:\, n\in\omega,\ a\text{ an assignment}\}.
\end{equation}

The following property is trivial, but important: Fix $n$. Then there is a
finite set $I$ such that we can determine the value that is assigned to $x_n$
by an assignment $a$ provided we know the values $(t_i\circ a)_{i\in I}$.  We
denote this by the following: For a canonical term sequence $\bar t$,
\begin{equation}\label{eq:determined}
  \text{each $x_n$ is determined by finitely many $t_i$.}
\end{equation}
(Proof: Let $l$ be the maximum of the heights of the nodes in $F_n$. Set
$I=\{0,\dots,l\}$.)

In the example above, $x_0$ is determined by $t_0$, and $x_1$ by
$(t_0,t_1,t_2)$, but not by $(t_0,t_1)$.

Let $T'\subseteq T$ be a perfect subtree, 
and call the canonical term sequence $(t'_0,t'_1,\dots)$,
written as terms in the variables $x'_0,x'_1,\dots$.
In the example  of Figure~\ref{fig:tree}, we get:
\[
  t'_0=1,
  \quad
  t'_1=1,
  \quad
  t'_2=x'_0.
\]
The fronts $F'_n$ ``refine'' $F_n$: 
If $t\in F'_n$, then $t\geq s$ for a unique $s\in F_n$.
So the variables $(x'_0,\dots,x'_l)$ give at least as much information (about
the branch) as $(x_0,\dots,x_l)$.  In other words, we can calculate the value
of $x_n$ given the values $(x'_0,\dots x'_n)$, and we write this dependence as
a term $\phi_n(x'_0,\dots x'_n)$.  This defines a function (or: term sequence)
$\phi$  that assigns to each variable $x_n$ a term $\phi_n(x'_0,\dots,x'_n)$.
We will call $\phi$ a substitution. So for every assignment $a$ of the
variables $x'$, we get the same result when we apply $a$ to the term sequence
$\bar t'$ as we get when we apply $\phi\circ a$ to $\bar t$.
In other notation, $\bar t'=\bar t\circ \phi$.

In the example, the substitution $\phi$ has the following values:
\[
  x_0=\phi_0(x'_0)=1,\quad x_1=\phi_1(x'_0,x'_1)=x'_0, \quad \dots
\]
It is easy to check that, e.g., $t'_2(x'_0, x'_1,x'_2)=x'_0$ is indeed the same
as $t_2(x_0,x_1,x_2)$ after applying the substitution $\phi$,
i.e., $t_2'=t_2\circ \phi$: 
\[
  t_2
  =
  \begin{cases}
         x_2&\text{if }x_0=0\text{ and }x_1=0\\
         1  &\text{if }x_0=0\text{ and }x_1=1\\
         x_1&\text{otherwise.}
  \end{cases}
  \quad\quad
  t_2\circ \phi=
  \left\{\begin{aligned}
         \phi_2&\text{ if }1=0\text{ and }x'_0=0\\
         1  &\text{ if }1=0\text{ and }x'_0=1\\
         x'_0&\text{ otherwise}
  \end{aligned}\right\}
  =x'_0.
\]

Also, each $x'_j$ is determined by finitely many $\phi_i$.
This means: For each $j$ there is a finite set $I$ such that the following
holds: If $a,b$ are assignments of $\{x'_0,x'_1,\dots\}$ that map $x'_j$ to
different values, then
$(\phi_i\circ a)_{i\in I}\neq (\phi_i\circ b)_{i\in I}$.
(Proof: Pick $l\in\omega$ such that each node in $F_l$ is longer than every
node in $F'_{j+1}$, and set $I=\{0,\dots,l\}$.)

So far, we used different variable symbols ($x_i$ and $x_i'$) for variables
used in $\bar t$ and $\bar t'$ (in the hope to make the concept of substitution
a bit clearer). Of course this is not necessary, and we will only
use $x_i$ in the following.
We will see that the following partial order $S^*$ is equivalent to Sacks
forcing: $S^*$ consists of sequences of terms $(t_i)_{i\in\omega}$
using the variables $x_j$ ($j\in\omega$) such that
\begin{eqpar}\label{Sfirst}
  \begin{enumerate}
    \item $t_i$ depends only on $x_j$ with $j\leq i$, and
    \item each $x_j$ is determined by finitely many $t_{i}$.
  \end{enumerate}
\end{eqpar}
The order is defined as follows:  $\bar t'$ is stronger than $\bar
t$, if there is a substitution $\phi$ such that $\bar t'=\bar t\circ \phi$ and
\begin{eqpar}\label{orderfirst}
  \begin{enumerate}
    \item $\phi_i$ only depends on $x_j$ with $j\leq i$, and
    \item each $x_j$ is determined by finitely many $\phi_i$.
  \end{enumerate}
\end{eqpar}

It is easy to check that $\leq$ is reflexive and transitive; and that $\circ$
is associative: The identity substitution witnesses $\bar t\leq \bar t$; and if
$\bar t'=\bar t\circ \phi'$ and $\bar t''=\bar t'\circ \phi'$ then $\bar
t''=(\bar t\circ \phi)\circ \phi'=\bar t\circ (\phi\circ \phi')$.

We could omit \eqref{orderfirst}(ii): If $\phi$ is any substitution, and if
$\bar t$ and $\bar t\circ \phi$ both are in  $S^*$, then $\phi$ satisfies~(ii)
anyway.

Subsitutions (as defined in \eqref{orderfirst}) are obviously exactly the same
as conditions in $S^*$ (as defines in \eqref{Sfirst}). This fact is not deep or
of any real importance, but it will simplify our notation. So let us describe
this effect once more:

Assume that $\bar s$ and $\bar t$ both are conditions in $S^*$. We can
interpret $\bar t$ as a substitution $\phi$ such that $\phi_n=t_n$. (I.e.,
$t_n(x'_0,\dots,x'_n)$ calculates the value of $x_n$.)   
Then $\bar s\circ \bar t$ is
again element of $S^*$ (and stronger than
$\bar s$). On the other hand, if $\bar t'$ is stronger than $\bar s$,
then this is witnessed by a substitution $\phi$, which we can in
turn interpret as element of $S^*$.

We can interpret a $\bar t\in S^*$ as continuous function from $2^\omega$ to
$2^\omega$, and map $\bar t$ to its image, or to the associated tree:

\begin{Lem}\label{lem:guri}
  Let $\Psi$ map $\bar t\in S^*$ to 
  $\{(\bar t\circ a)\restriction n:\, n\in\omega,a\text{ an assignment}\}$.
  Then $\Psi$ is a surjective complete embedding (in particular order
  preserving) from $S^*$ into Sacks forcing.
\end{Lem}

\begin{proof}
  $\Psi(\bar t)$ is a perfect tree: Pick any
  $s=(\bar t\circ a)\restriction n\in \Psi(\bar t)$.
  Note that $t_0,\dots t_{n-1}$ use a finite set $A$ of
  variables. Pick  $x_j\notin A$. Then $x_j$ is determined
  by $t_0,\dots t_{l-1}$ for some $l$. Pick assignments $b,c$
  extending $a\restriction A$ such that
  $x_j\circ b\neq x_j\circ c$.
  Then $(\bar t\circ b)\restriction l\neq (\bar t\circ c)\restriction l$
  (otherwise they would determine the same value for $x_j$), so we get two
  incomparable nodes in $\Psi(\bar t)$ both extending $s$.

  We see from \eqref{eq:blub2} that
  $\Psi$ is surjective. It is clear that $\Psi$ preserves $\leq$.

  $\Psi$ preserves $\incomp$:
  Assume that
  $\Psi(\bar t)$ and $\Psi(\bar s)$ both contain the perfect tree $T$.
  By thinning out $T$, we can assume the following: If $l$ is the  
  length of a node in $F_n$, then $t_1(\bar x),\dots,t_l(\bar x)$
  determine $x_n$, and the same holds for $\bar s$.
  Let $\bar r(\bar x')\in S^*$ be the canonical sequence of $T$.
  So $x'_0\dots x'_{n-1}$ determine a node in $F_n$, and therefore
  sufficiently many $t_1,\dots,t_{l-1}$ to determine $x_n$.
  This defines a substitution $\phi$ witnessing that $\bar r$
  is stronger than $\bar t$. The same applies to $\bar s$.
\end{proof}

Of course $\Psi$ is not injective. For example, if we simply interchange $x_0$
and $x_1$ in a suitable sequence $\bar t$, then we can still get a valid term
sequence (different from the original one), but the image under $\Psi$ will be
the same.

\begin{figure}[tb]
    \begin{center}
      \scalebox{0.8}{\input{order.pstex_t}}
    \end{center}
  \caption{\label{fig:order}%
  We use a canonical ordering of $\omega\times \omega$.
  The node $(1,2)$ corresponds to $7$. The nodes $(n,m)$ smaller than 
  $(1,2)$ all satisfy $n+m\leq 1+2<4$, as in \eqref{eq:blug}.
  }
\end{figure}
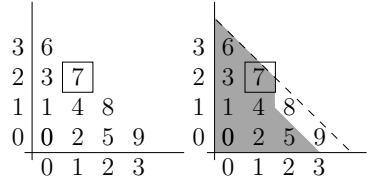
In $S^*$, the index set of the term sequence is $\omega$. We will later need
$\omega\times \omega$-sequences; so we will just identify $\omega$ with
$\omega\times\omega$, using a canonical order. See Figure~\ref{fig:order}.  
\begin{equation}
  \tau:\omega\times\omega\to\omega\text{ is defined by }
  \tau(n,m)=n+\frac{1}{2}(n+m)(n+m+1).
\end{equation}
The bijection $\tau$ 
defines a linear order of $\omega\times\omega$ of order type $\omega$:
\begin{equation}
  (i,j)\unlhd (n,m)\text{ iff }\tau(i,j)\leq \tau(n,m)
\end{equation}
We will later use the following trivial fact:
\begin{equation}\label{eq:blug}
  \text{If }i+j<n\text{ and }(i',j')\unlhd (i,j)\text{ then }i'+j'<n.
\end{equation}

We now rewrite $S^*$ in the form of $\omega\times \omega$-sequences:
\begin{Def}\label{def:Qs}
  $Q_*$ consists of squares of terms $(t_{n,m})_{n,m\in\omega}$ 
  using the variables $x_{i,j}$ ($i,j\in\omega$) such that
  \begin{enumerate}
    \item $t_{n,m}$ depends only on $x_{i,j}$ with $(i,j)\unlhd (n,m)$, and
    \item each $x_{i,j}$ is determined by finitely many $t_{n,m}$.
  \end{enumerate}
  The order is defined as follows:  $\bar t$ is stronger than $\bar
  s$, if there is a condition $\phi\in Q_*$ such that
  $\bar t=\bar s\circ \phi$.
\end{Def}

Since $Q_*$ is isomorphic to $S^*$, Lemma~\ref{lem:guri} gives us:
\begin{Cor}
  $Q_*$ is forcing equivalent to Sacks forcing.
\end{Cor}

We now add the a formal definition of term, assignment and substitution:
\begin{Def}\label{def:terms}
  \begin{itemize}
    \item Let $X$ be a set. We will call an element $v\in X$
      a variable (or: variable symbol).
      We will interpret $v$ as a ``binary variable'',
      i.e., the value of $v$ is 0 or 1.
      In this paper, we will use $X=\{x_i:\, i\in\omega\}$
      and $X=\{x_{i,j}:\, i,j\in\omega\}$.
    \item An $X$-term $t$ consists of\footnote{Formally we could
      let $t$ be a triple $(X,(v_0,\dots, v_{l-1}),f)$,
      to guarantee that $X$ is disjoint to the terms built from it.}
      a sequence $(v_0,\dots, v_{l-1})$
      for some $0\leq l<\omega$ and $v_i\in X$, together with a 
      function $f:2^l\to 2$. (So for $l=0$, the sequence of variables
      is empty and the term is a constant.) 
      We usually write terms as $t(v_0,\dots, v_{l-1})$.
      Abusing notation, we identify the variable $v$ with the
      ``identity term'' corresponding to $(v),\textrm{Id}$.
    \item
      An assignment $a$ is a function $X\to 2$. Assignments extend to
      all $X$-terms in the natural way. In other words,
      given an assignment $a$, we can apply $a$ to a term $t$ 
      to get an element of $2$. We denote the result of applying 
      $a$ to a term (or variable) $t$ by $t\circ a\in 2$.
    \item
      Similarly, a substitution $\phi$ maps $X$ to $X$-terms.
      Equivalently, a substitution is a sequence $(\phi_v)_{v\in X}$
      of $X$-terms indexed by $X$.
      Again, we can extend a substitution to act on all $X$-terms,
      and we write $t\circ \phi$ for the result.
      We can also apply substitutions to sequences $\bar t=(t_v)_{v\in X}$
      of terms (indexed by $X$), the result $\bar t \circ \phi $ 
      is another sequences of terms indexed by $X$.
      The application of substitutions is associative:
      For term sequences $\bar r$, $\bar s$ and 
      $\bar t$, all indexed by $X$, we get 
      $ \bar r\circ (\bar s\circ \bar t)=(\bar r\circ \bar s)\circ \bar t$.
    \item The variable (or term) $s$ ``is determined by the terms $t_0,\dots t_n$''
      means that 
      \[
        (t_0\circ a,\dots t_n\circ a)=(t_0\circ b,\dots ,t_n\circ b)
        \text{ implies }
        s\circ a=s\circ b
      \]
      for all assignments $a$ and $b$. In other words, if we know the 
      value of $t_0,\dots t_n$, we can infer the value of $s$.
    \item According to our formal definition, two terms 
      that depend on different variables
      are distinct 
      (even if these variables are not relevant). However, we will only be
      interested in terms ``as functions'', i.e., modulo
      the following equivalence relation: 
      $t=^* s$ means that $t\circ a=s\circ a$ for all assignments $a$.
      In particular, the last $=$ sign in Definition~\ref{def:Qs} really means
      $=^*$ etc.
  \end{itemize}
\end{Def}

\section{A simple case}\label{sec:easy}

In the rest of the paper, $\delta$ always denotes a countable limit ordinal.

In this section, we construct a proper, NNR countable support iteration and
argue that the limit adds a real that it is similar to a Sacks real (i.e., 
it adds a
generic object for a forcing that looks in some way similar to the 
$Q_*$ defined in the previous section).
In the rest of the paper, we deal with an analog (but notationally
more complicated) construction that actually adds a Sacks real.

So the purpose of this section is to give some idea of the constructions we use
to prove Theorem~\ref{thm:main}, using a somewhat simplified notation.  The
reader who does not feel the need of such an introduction can safely continue
with the next section. 

We do not give any proofs in this section, but refer to the proofs of the more
general statements. {\bf Caution:} We use the same symbols for the simpler
objects in this section and for the analog constructions in the rest of the
paper.

The forcing iteration will start with a preparatory forcing $\tilde P$,
followed by $\n Q_0,\n Q_1, \dots$.  $\tilde P\ast P_n$ stands for $\tilde
P\ast Q_0\ast \dots\ast Q_{n-1}$. We will also use the countable support limit
of $\tilde P\ast P_n$.  Since all forcings are proper, this countable support
limit is the same as $\tilde P\ast P_\omega$, where $P_\omega$ is the
$\tilde P$-name for the countable support limit of the $P_n$.

The preparatory forcing adds
cofinal 
subsets $\nu_{\delta,n,m}\subseteq \delta$ of of order type $\omega$
for every limit ordinal $\delta<\om1$ and $n,m\in\omega$.
In more detail:
\begin{Def}
  A condition $\tilde p$ in $\tilde P$ consists of a 
  limit ordinal $\height(\tilde p)\in\om1$ and a sequence 
  $(\nu_{\delta,n,m})_{0<\delta<\height(\tilde p), n,m\in\omega}$,
  such that 
  $\nu_{\delta,n,m}\subseteq \delta$ is cofinal and has order type $\omega$, and
  $\nu_{\delta,n,m_1}$ and $\nu_{\delta,n,m_2}$ are disjoint
  for $m_1\neq m_2$.
  $\tilde P$ is ordered by extension.
\end{Def}

So $\tilde P$ is $\sigma$-closed.

\begin{Def}
  $Q_0$ is (the $\tilde P$-name) for $2^{<\om1}$, ordered by extension.
\end{Def}
So $Q_0$ is $\sigma$-closed as well, and adds the generic sequence
$\eta_0\in 2^{\om1}$.

Given $\tilde P\ast P_n=\tilde P\ast Q_0\ast \dots\ast Q_{n-1}$ such that
$Q_{n-1}$ adds the generic sequence $\eta_{n-1}\in2^{\om1}$, we define the
$\tilde P\ast P_n$-name $Q_n$ (see also Figure~\ref{fig:simple}(a)):
\begin{Def}
  Let $q$ be a partial function from $\om1$ to $2$, 
  $\delta\subseteq \dom(q)$. $q$ and $\eta_{n-1}$ cohere
  at $\delta+m$, if $\eta_{n-1}(\delta+m)=q(\alpha)$
  for all but finitely many $\alpha\in \nu_{\delta,n-1,m}$.
  Abusing notation, we just say $q$ coheres with $\eta_{n-1}(\delta+m)$.
  \\
  We set $q\in Q_n$, if $q\in 2^{<\om1}$ and
  $q$ coheres with $\eta_{n-1}(\delta+m)$
  for all $\delta\leq \dom(q),m\in\omega$.
\end{Def}
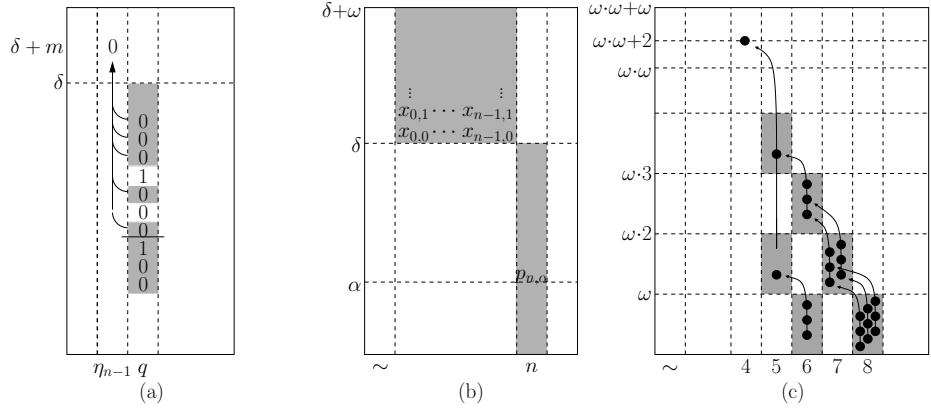
\begin{figure}[tb]
      \newcommand{\mylabelsimplecoherence}{(a)}
      \newcommand{\mylabelR}{(b)}
      \newcommand{\mylabelsimpleexample}{(c)}
  \hfill
  \begin{minipage}[t]{.3\textwidth}
    \begin{center}
      \scalebox{0.4}{\input{simplecoherence.pstex_t}}
    \end{center}
  \end{minipage}
  \hfill
  \begin{minipage}[t]{.3\textwidth}
    \begin{center}
      \scalebox{0.4}{\input{R.pstex_t}}
    \end{center}
  \end{minipage}
  \hfill
  \begin{minipage}[t]{.3\textwidth}
    \begin{center}
      \scalebox{0.4}{\input{simpleexample.pstex_t}}
    \end{center}
  \end{minipage}
  \hfill
  \caption{\label{fig:simple}
     \mylabelsimplecoherence{}~$q$ coheres with $\eta_{n-1}(\delta+m)$.
     The gray area indicates $\nu_{\delta,n-1,m}$.
     \mylabelR{}~An element of $R$: 
     $p_{n,\delta+m}=x_{n,m}$. The term $p_{n,\alpha}$ only
     depends on $x_{i,j}$ with $i<n$.
     \mylabelsimpleexample{}~Filling in the term $x_{4,2}$
     at at various positions: In the bottom row,
     $I_{4,2}$ contains 6 and 8; the $1\times \omega$-blocks 
     where the terms $x_{4,2}$ are added (indicated by the gray area)
     propagates up-left.}
\end{figure}

\begin{Lem}\label{lem:jkwru}
  The following is forced by $P_n$:
  \begin{enumerate}
    \item If $q\in Q_n$, $q'\in 2^{\dom(q)}$, and
      $q'(\alpha)=q(\alpha)$ for
      all but finitely many $\alpha\in\dom(q)$, then $q'\in Q_n$.
    \item $Q_n$ is separative,%
      \footnote{That is, for every $p\in Q$ there are $q_1,q_2\leq p$ such that
        $q_1\incomp q_2$.}
       and adds a generic sequence
       $\eta_n\in 2^{\om1}$ defined by $\bigcup_{q\in G(n)}q$.
    \item If $n>0$, then $Q_n$ is not $\sigma$-closed.
    \item 
      $Q_n$ is proper and NNR, i.e., $Q_n$ adds no new real. 
  \end{enumerate}
\end{Lem}
For a proof, see Lemmas~\ref{lem:basicQ},~\ref{lem:newreal}
and~\ref{lem:proper}. Note that (i)--(iii)  are very easy, and (iv) is
straightforward (but a bit cumbersome).

\begin{Lem}
  $\tilde P\ast P_\omega$ adds a new real. In particular,
  the $\tilde P$-generic element $\bar \nu$ together with 
  $(\eta_n(m))_{n,m\in\omega}$ determines the generic filter.
\end{Lem}
\begin{proof}
  If we know $\nu_{\omega,n-1,m}$ and
  $\eta_n(l)$ for all $l\in\omega$, then we can determine
  $\eta_{n-1}(\omega+m)$.
  So if we know all $\nu_{\delta,n,m}$ and all $\eta_n(m)$
  for $\delta<\om1,n,m\in\omega$,
  then we can by induction on $\alpha<\om1$ calculate all $\eta_n(\alpha)$
  for $n\in\omega$.
\end{proof}

We now define a dense subforcing $R$ of $\tilde P\ast P_\omega$. See
Figure~\ref{fig:simple}(b). We use the notion of variable, term, assignment and
substitution, as in Definition~\ref{def:terms}, for the set of variables
$X=\{x_{i,j}:\, i,j\in\omega\}$.

\begin{Def}
  $R=\bigcup_{\delta<\om1}R_{\delta+\omega}$.
  A condition $p$ in $R_{\delta+\omega}$ consists of $\tilde p$ and $\bar p$ such that
  \begin{itemize}
    \item $\tilde p\in \tilde P$, $\height(\tilde p)=\delta+\omega$.
    \item $\bar p=(p_{n,\alpha})_{n\in\omega,\alpha\in\delta+\omega}$.
    \item $p_{n,\delta+m}$ is the term $x_{n,m}$. 
    \item For $\alpha<\delta$,
      $p_{n,\alpha}$ is a term using only variables $x_{i,j}$ with $i<n$.
    \item For $\alpha\leq \delta$ limit and $n,m\in\omega$,
      $p_{n,\alpha+m}=p_{n+1,\zeta}$ for all but finitely
      many $\zeta\in\nu_{\alpha,n,m}$.
  \end{itemize}
\end{Def}
We identify two conditions $p$ and $q$ if $\tilde p=\tilde q$ and
$p_{n,\alpha}=^*q_{n,\alpha}$ for all $n\in\omega$, $\alpha<\delta$.


We can interpret $p\in R$ as a condition $(\tilde
p,p(0),p(1),\dots)$ in $\tilde P\ast P_\omega$: After forcing with $\tilde
P\ast P_n$, we have the generic sequences $(\eta_i)_{i<n}$. This defines a
canonical assignment of $x_{i,j}$ for $i<n$, namely $x_{i,j}:=\eta_i(\delta+j)$. This assignment evaluates $(p_{n,\alpha})_{\alpha<\delta}$
to a condition in $Q_n$ (assuming that $\tilde p$ is 
element of the $\tilde P$-generic filter), 
and we define $p(n)$ to be that condition.
Using this identification, we get:

\begin{Lem}
    $R$ is a dense subset of $\tilde P\ast P_\omega$.
\end{Lem}
For a proof, see Lemma~\ref{lem:dense}.
The proof is again a bit cumbersome, and uses similar arguments (chains of
countable elementary submodels) as the proof of \ref{lem:jkwru}(iv).

Note the following simple properties for $p\in R_{\delta+\omega}$:
\begin{itemize}
  \item If $\delta=\omega$, we get
      $(\forall n,m)(\exists^\infty k) p_{n+1,k}=x_{n,m}$.
  \item If  $\delta=\omega+\omega$, we get
      $(\forall n,m)(\exists^\infty k) p_{n+2,k}=x_{n,m}$.
  \item If $\delta=\omega\cdot\omega$, then we get
      $(\forall n,m)(\exists^\infty n') (\exists^\infty m') p_{n',m'}=x_{n,m}$.
\end{itemize}
Actually, the last item holds for all $\delta\geq \omega\cdot\omega$, which can
easily be seen by induction; and we get some kind of converse as well:
\begin{Lem}\label{lem:3515123}
  Assume that $(r_{i,j})_{i,j\in\omega}$ is a matrix of terms
  such that
  \begin{enumerate}
    \item $r_{i,j}$ depends only on $x_{n,m}$ with $n< i$,
    \item $(\forall n,m)(\exists^\infty n') (\exists^\infty m') r_{n',m'}=x_{n,m}$
    \item $(\forall n)(\exists^\infty k) r_{n,k}=0$.
  \end{enumerate}
  Then there is a $p\in R_{\omega\cdot\omega+\omega}$ such that
  $r_{i,j}=p_{i,j}$ for all $i,j\in\omega$.
\end{Lem}
\begin{proof}[Proof (sketch)]
  We have to define a suitable $\tilde p$ 
  (i.e., the sequence $(\nu_{\alpha,n,m})_{\alpha\leq \omega,n,m\in\omega}$)
  as well as $p_{i,\alpha}$ for
  $i\in\omega$ and $\omega\leq \alpha< \omega\cdot\omega$.

  We deal with one variable after the other,
  see Figure~\ref{fig:simple}(c).
  Assume we are dealing with $x_{n,m}$.
  Set 
  \[
    I_{n,m}=\{n'\geq n+2:\, (\exists^\infty k) r_{n',k}=x_{n,m}\}.
  \]
  According to~(ii), $I_{n,m}$ is infinite.
  We define $\nu_{\omega\cdot\omega,n,m}\subseteq [\omega,\omega\cdot\omega[$
  such that
  \[
    \nu_{\omega\cdot\omega,n,m}\cap [i\cdot \omega,(i+1)\cdot \omega[
  \]
  contains a single element (not used so far) 
  if $i+m+1\in I_{n,m}$ and is empty otherwise.
  We set $r_{n+1,\alpha}=x_{n,m}$ for all $\alpha\in \nu_{\omega\cdot\omega,n,m}$;
  and propagate the $x_{n,m}$ diagonally down.

  We repeat the same construction for all the other $x_{n,m}$, and then
  set all the remaining terms $r_{n,\alpha}=0$. To get coherence for 
  these points as well, we just define the remaining
  $\nu$'s in a way so that they only point to $r$'s  that are 0.
  At height $\omega$, we use~(iii) to do this,
  at other heights we just have to make sure to leave enough space when 
  choosing the elements of $\nu_{\delta,n,m}$.
\end{proof}

We now describe how to ``stack'' a condition on top of another one to get a
stronger condition.  See Figure~\ref{fig:R}(b) for a graphical illustration.
\begin{itemize}
  \item If we ``cut away the bottom part'' of a condition
    $q\in R_{\delta+\delta'+\omega}$ at height $\delta$, then we
    get a condition $q'\in R_{\delta'+\omega}$. Formally we can define $q'$ as 
    follows:
    \begin{itemize}
      \item $\beta\in \tilde q'_{\alpha,n,m}$ iff $\delta+\beta\in \tilde q_{\delta+\alpha,n,m}$.
      \item $q'_{n,\alpha}=q_{n,\delta+\alpha}$.
    \end{itemize}
    We denote this $q'$ by $q\restriction [\delta,\delta+\delta'+\omega]$.
  \item 
    We can stack any condition $q'\in R_{\delta'+\omega}$ 
    on top of some condition $p\in R_\delta+\omega$, resulting in some $q\in
    R_{\delta+\delta'+\omega}$ such that
    $q\restriction [\delta,\delta+\delta'+\omega]=q'$. Formally, $q$ is defined
    as follows:
    \begin{itemize}
      \item $\tilde q\restriction \delta=\tilde p$.
      \item For $\alpha<\delta'$,  we set 
        $\delta+\beta\in\tilde q_{\delta+\alpha,n,m}$
        iff $\beta\in\tilde q'_{\alpha,n,m}$.
      \item For $\alpha<\delta'$,  we set 
        $q_{n,\delta+\alpha}=q'_{n,\alpha}$.
      \item We define the substitution $\phi$ by
        $\phi_{n,m}=q'_{n,m}$, and set
        $q_{n,\alpha}=p_{n,\alpha}\circ \phi$ 
        for all $\alpha<\delta$.
    \end{itemize}
    We denote this $q$ by $p\Lsh q'$.
\end{itemize}
It is clear that $p\Lsh q'\leq p$ 
(interpreted as element of $\tilde P\ast P_\omega$).
The converse is true as well: 
\begin{eqpar}\label{eq:gui}
  If $q\leq p$, then either $q=p$ or
  $q=p\Lsh q'$, where $q'=q\restriction [\height(p),\height(q)]$.
\end{eqpar}

The proof of~\eqref{eq:gui} uses the following simple
fact: For every  $p\in R$,
\begin{eqpar}
  every partial assignment of every finite
  $A\subseteq \{x_{i,j}:i,j\in \omega\}$ is compatible
  with $p$.
\end{eqpar}
In other words, if $f:\omega\times\omega\to 2$ is a finite partial function,
then it is compatible with $p$ (interpreted as element of $\tilde P\ast
P_\omega$) that $\eta_n(\delta+m)=f(n,m)$ for all $(n,m)\in\dom(f)$.

Given a $p\in R_\delta$ (we assume $\delta\geq \omega\cdot \omega)$, 
we can map $p$ to the square of terms
$\sigma(p)=(p_{n,m})_{n,m\in\omega}$.
Then $\sigma$ maps $R$ to $Q_{**}$ in an order preserving way,
where $Q_{**}$ is defined as follows:

\begin{Def}
  $Q_{**}$ is the set of all sequences $\bar t=(t_{i,j})_{i,j\in\omega}$
  of terms such that
  \begin{enumerate}
    \item $t_{n,m}$ depends only on $x_{i,j}$ with $i<n$,
    \item $(\forall n,m)(\exists^\infty n')(\exists^\infty m') t_{n',m'}=x_{n,m}$,
  \end{enumerate}  
  $\bar t\leq \bar s$, if there is a substitution $\phi$ such that:
  \begin{enumerate}
    \setcounter{enumi}{2}
    \item $t_{i,j}=s_{i,j}\circ \phi$ for all $i,j\in\omega$.
    \item $\phi_{n,m}$ only depends on
      $x_{i,j}$ with $i\leq n$.
  \end{enumerate}
\end{Def}

\begin{Lem}
  $R$ (or equivalently: $\tilde P\ast \tilde P_\omega$) adds a
  generic filter for $Q_{**}$.
\end{Lem}
Note that $Q_{**}$ looks somewhat similar to the $Q_*$ defined in the previous
section. For  $Q_*$ instead of $Q_{**}$, the Theorem is the main part of
Theorem~\ref{thm:main}. In the rest of the paper, we will modify the
constructions so that we actually end up with  $Q_*$ instead of $Q_{**}$.

\begin{proof}[Proof (sketch)]
  We already mentioned that 
  $\sigma: R\to Q_{**}$ preserves $\leq$. Assume that
  $G$ is  $R$-generic over $V$, and define
  \[
    G_{**}=\{\bar t\in Q_{**}:\, (\exists p\in G)\, \sigma(p)\leq \bar t\}.
  \] 
  It is enough to show the following:
  \begin{eqpar}\label{eq:wqwrtw}
    For $p\in R$ there is a  $\bar s\leq_{Q_{**}} \sigma(p)$
    such that for all $\bar t\leq \bar s$  then there is 
    an $q\leq_R p$ such that $\sigma(q)\leq \bar t$.
  \end{eqpar}
  Then the Lemma follows:
  First note 
  that $ G_{**}$ does not
  contain incompatible elements, since $\sigma$ is order preserving.
  Now assume that $D\subseteq Q_{**}$ is dense, 
  and (towards a contradiction) that $p$ forces that $G_{**}$ does not 
  meet $D$. Then pick some $\bar t\leq \bar s$ in $D$ and some 
  $q$ as above, contradiction.

  To show~\eqref{eq:wqwrtw}, we define  $\bar s\in Q_{**}$
  via the substitution $\phi$ witnessing $\bar s\leq \sigma(p)$,
  defined as follows:
  For each $n$, let $(\phi_{n,m})_{m\in\omega}$ enumerate (with infinite
  repetitions) the constant term $0$ and all variables
  $x_{i,j}$ with $i<n$. So $\phi$ maps $x_{n,m}$ to 
  the term $\phi_{n,m}$.

  Now pick any $\bar t$ that is stronger than $\bar s$, 
  witnessed by some substitution $\psi$. 
  Note that  $\phi\circ \psi$ satisfies the requirements of 
  Lemma~\ref{lem:3515123}. So there is a 
  $q'\in R$ such that $\sigma(q')=\phi\circ \psi$. Then
  $p\Lsh q'$ is as required.
\end{proof}

In the rest of the paper, we will modify the constructions of this section in
such a way that we end up with $Q_*$ instead of $Q_{**}$. It turns out that
this does not require any new concepts, just a more awkward notation.

\section{The NNR iteration}

In the rest of the paper, $\delta$ always denotes a countable limit ordinal.

First we define a $\sigma$-closed preparatory forcing $\tilde P$, which gives
us for every limit
$\alpha\in\om1$ a subset of $\alpha$ of order type $\omega$ and
some simple coding sequences.

\begin{Def}
  $\tilde p\in \tilde P$ if for some $\height(\tilde p)<\om1$,
  $\tilde p$ consists of sequences
  \[
    \nu_{\delta,n,m},\ j_{\delta,n,m},\text{ and }f_{\delta,n,m,k}
    \text{ for } \delta<\height(\tilde p),\ m,n,k\in\omega,
  \]
  such that
  \begin{itemize}
    \item each $\nu_{\delta,n,m}$ is a cofinal, unbounded subset of $\delta$
      of order type $\omega$.
    \item $m_1\neq m_2$ implies that
      $\nu_{\delta,n,m_1}$ and $\nu_{\delta,n,m_2}$ are disjoint.
    \item $j_{\delta,n,m}$ is an increasing function from $\omega$ to $\omega$.
    \item $f_{\delta,n,m,k}$ is a surjective function from 
       $2^{[j_{\delta,n,m}(k),j_{\delta,n,m}(k+1)-1]}$ to $2$.
  \end{itemize}
  $\tilde P$ is ordered by extension. 
\end{Def}

\begin{Lem}
  $\tilde P$ is $\sigma$-closed, and forces that $2^{\al0}=\al1$.
\end{Lem}

\begin{proof}
  In the $\tilde P$-extension $V'$
  define $A_\delta$ by 
  $l\in A_\delta\text{ iff } \delta+2l\in \nu_{\delta+\omega,0,0}$.
  By a simple density argument, 
  $\{A_\delta:\, \delta<\om1\}$  contains all old reals and
  therefore all reals.
\end{proof}

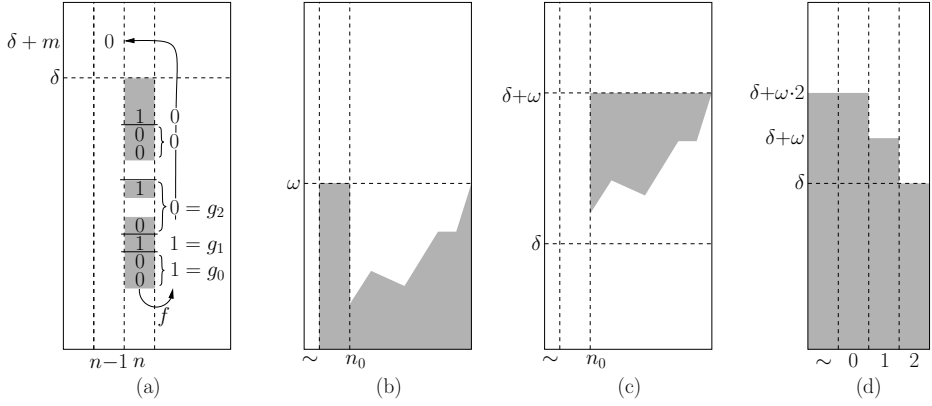
\begin{figure}[tb]
      \newcommand{\mylabelcoherence}{(a)}
      \newcommand{\mylabelpossible}{(b)}
      \newcommand{\mylabelimpossible}{(c)}
      \newcommand{\mylabeltypicalfinite}{(d)}
  \hfill
  \begin{minipage}[t]{.24\textwidth}
    \begin{center}
      \scalebox{0.4}{\input{coherence.pstex_t}}
    \end{center}
  \end{minipage}
  \hfill
  \begin{minipage}[t]{.24\textwidth}
    \begin{center}
      \scalebox{0.4}{\input{possible.pstex_t}}
    \end{center}
  \end{minipage}
  \hfill
  \begin{minipage}[t]{.24\textwidth}
    \begin{center}
      \scalebox{0.4}{\input{impossible.pstex_t}}
    \end{center}
  \end{minipage}
  \hfill
  \begin{minipage}[t]{.24\textwidth}
    \begin{center}
      \scalebox{0.4}{\input{typicalfinite.pstex_t}}
    \end{center}
  \end{minipage}
  \hfill
  \caption{\label{fig:P} \mylabelcoherence{}~$\eta_{n}$ coheres with
     $\eta_{n-1}(\delta+m)$ above $k_0=2$.
     $\nu_{\delta,n-1,m}$ is indicated by the gray
     area, $j_{\delta,n-1,m}$ corresponds to the partition of this
     area, $f_{\delta,n-1,m,k}$ calculates $g_{\delta,n-1,m,k}$.
     \mylabelpossible{}~There are conditions which determine the gray values,
     but these conditions are not dense. 
     \mylabelimpossible{}~No condition can determine all gray values.
     \mylabeltypicalfinite{}~A typical condition in $\tilde P\ast P_3$:
     The gray area indicates the domain,  all the values are 
     ``constants''.}
\end{figure}

Fix $\delta,n,m,k$.
Given  $\nu_{\delta,n-1,m}$, $j_{\delta,n-1,m}$
and $f_{\delta,n-1,m,k}$,
we define the function
$g_{\delta,n-1,m,k}: 2^\delta\to \{0,1\}$ the following way,
cf.~Figure~\ref{fig:P}(a): 

Fix $\eta_{n}\in 2^\delta$.
For $i\in\omega$, let $\zeta_i$ be the
$i$-th element of $\nu_{\delta,n-1,m}$. Set $b_i=\eta_{n}(\zeta_i)$.
So $\bar b=(b_i)_{i\in\omega}\in 2^\omega$. 
Look at $\bar b\restriction [j_{\delta,n-1,m}(k),j_{\delta,n-1,m}(k+1)-1]$.
This is a $0$-$1$-sequence of appropriate length, so we can
apply $f_{\delta,n-1,m,k}$. We call the result $g_{\delta,n-1,m,k}(\eta_{n})$.
To summarize:
Let $\zeta_i$ be the $i$-th element of $\nu_{\delta,n-1,m}$. Then
we define
\[
  g_{\delta,n-1,m,k}(\eta)=f_{\delta,n-1,m,k}
  ((\eta(\zeta_i))_{j_{\delta,n-1,m}(k)\leq i <j_{\delta,n-1,m}(k+1)}).
\] 

We will be interested in sequences $(\eta_n)_{n\in\omega}$ that cohere
with respect to $g_{\delta,n,m,k}$.
We again refer to~Figure~\ref{fig:P}(a):
\begin{Def}\label{def:cohere}
  Let $\eta_{n-1}$ and $\eta_{n}$ be partial functions from
  $\omega_1$ to 2, $\delta+m\in \dom(\eta_{n-1})$,
  $\delta\subseteq\dom(\eta_{n})$, $k_0\in\omega$. We say that
  $\eta_{n-1}$ and $\eta_{n}$ cohere at $\delta+m$ above $k_0$, if
  $\eta_{n-1}(\delta+m)=g_{\delta,n-1,m,k}(\eta_{n})$
  for all  $k\geq k_0$. We say that
  $\eta_{n-1}$ and $\eta_{n}$ cohere at $\delta+m$, if they
  cohere above some $k_0$. 
  Abusing notation, we also say that $\eta_{n}$ coheres with
  $\eta_{n-1}(\delta+m)$.
\end{Def}

Let $\tilde G$ be $\tilde P$-generic over $V$, and define
in $V[\tilde G]$ the forcing notion $Q_0$:

\begin{Def}
  $p\in Q_0$ iff $p\in 2^{\height(p)}$ for some $\height(p)<\om1$.
  $Q_0$ is ordered by extension.
\end{Def}

So $Q_0$ is $\sigma$-closed and adds the generic $\eta_{0}\in
2^{\om1}$.
Assume that $n\geq 1$,
$\tilde P\ast P_n=\tilde P\ast Q_0\ast \dots \ast Q_{n-1}$, and 
$Q_{n-1}$ adds the generic sequence $\eta_{n-1}\in 2^{\om1}$.
Let $\tilde G\ast G_n$ be $\tilde P\ast P_n$-generic over
$V$. In $V[\tilde G\ast G_n]$, we define $Q_n$ the following way:

\begin{Def}
  $p\in Q_n$ iff $p\in 2^{\height(p)}$ 
  for some limit ordinal $\height(p)<\om1$ and
  $\eta_{n-1}$ and $p$ cohere everywhere, i.e.
  at all $\delta+m<\height(p)+\omega$.
  $Q_n$ is ordered by extension.
\end{Def}

\begin{notation}
  We denote the $\tilde P$-generic filter by $\tilde G$, we set
  $\tilde P*P_n=\tilde P\ast Q_0\ast \dots \ast Q_{n-1}$,
  with generic filter
  $\tilde G\ast G_n$. Since $\tilde P$ is proper,
  the countable support limit of $\tilde P\ast P_n$ is the same as
  $\tilde P*P_\omega$, where $P_\omega$ is the $\tilde P$-name
  for the countable support limit of the $P_n$. 
  The generic filter of $\tilde P*P_\omega$ is 
  denoted by $\tilde G\ast G_\omega$;
  and $G(n)$ is the $Q_n$-generic filter (a $\tilde P*P_{n+1}$-name,
  or equivalently a $\tilde P*P_{n}$-name for a $Q_n$-name).
\end{notation}

\begin{Lem}\label{lem:basicQ} The following is forced by $\tilde P\ast P_n$:
  \begin{enumerate}
    \item Conditions can be finitely modified: If $q\in Q_n$, $q'\in
      2^{\dom(q)}$ and
      $q'(\alpha)=q(\alpha)$ for
      all but finitely many $\alpha\in\dom(q)$, then $q'\in Q_n$.
    \item
      If $p\in Q_n$ and $\beta>\height(p)$ is a limit ordinal, then there
      is a $q\leq p$ with $\height(q)=\beta$. In particular,
      $Q_n$ adds the generic object $\eta_n=\bigcup G(n)\in 2^{\om1}$
      (which in turn determines the generic filter $G(n)$).
    \item
      $Q_n$ is separative (and in particular nontrivial), and not
      $\sigma$-closed for $n\geq 1$.
  \end{enumerate}
\end{Lem}

\begin{proof}
  (i) is trivial.

  (ii)
  Let $(\delta_i,m_i)_{i\in\omega}$
  enumerate all pairs $(\delta,m)$ such that 
  $\height(p)<\delta\leq \beta$ and $m\in\omega$.
  Define an increasing sequence $p_i$ of partial functions
  from $\beta$ to $\{0,1\}$:

  Set $p_0=p$. For $i>0$, assume that
  $\dom(p_{i-1})=\height(p)\cup \bigcup_{j<i} \nu_{\delta_j,n-1,m_j}$.
  Then $\nu_{\delta_i,n-1,m_i}\cap \dom(p_{i-1})$ is finite:
  \begin{itemize}
    \item
      If $j<i$ and
      $\delta_j=\delta_i$, then $m_j\neq m_i$ and
      $\nu_{\delta_i,n-1,m_i}$ and $\nu_{\delta_j,n-1,m_j}$ are disjoint.
    \item
      If $j<i$ and $\delta_j\neq \delta_i$, then $\nu_{\delta_i,n-1,m_i}\cap 
      \nu_{\delta_j,n-1,m_j}$ are finite (since $\nu_{\delta,n-1,m}$ is
      a cofinal subset of $\delta$ of order type $\omega$).
    \item For the same reason, $\nu_{\delta_i,n-1,m_i}\cap\height(p)$ is finite.
  \end{itemize}
  Therefore we can extend $p_{i-1}$ to some $p_{i}$ by adding values at
  $\nu_{\delta_i,n,m_i}\setminus \dom(p_{i})$
  that cohere with $\eta_{n-1}(\delta_i+m_i)$.
  (Recall that $f_{\delta_i,n-1,m_i,k}$ is onto.)
  Set $p_\omega=\bigcup p_n$, and fill in arbitrary values (e.g., 0) at
  $\beta\setminus \dom(p_\omega)$. This gives a $q\leq p$ with $\height(q)=\beta$.

  (iii) follows from (i) and (ii).
\end{proof}

\begin{remark}
  For this proof, as well as for most of the following, the preparatory
  forcing $\tilde P$ is not necessary: The definition of $Q_n$ works for any
  reasonably defined sequences $\bar \nu$, $\bar j$, $\bar f$. Only in 
  Section~\ref{sec:notagain}
  we need that these sequences are generic. 
  (Guessing with, e.g., a $\diamondsuit$-sequence is not enough,
  as discussed in Section~\ref{sec:quotient}.)
\end{remark}

\begin{Lem}\label{lem:newreal}
  $(\eta_n)_{n\in\omega}$ is determined by
  $\tilde G$ and
  $(\eta_n(m))_{n,m\in\omega}$.
  In particular, $\tilde P\ast P_\omega$ adds a new real.
\end{Lem}
Of course, we do not use any particular property of the
countable support limit here. More generally, we get:

Assume $V'$ is an extension of $V$ that contains some $\tilde G$
and a sequence
$(G(n))_{n\in\omega}$ such that 
$\tilde G$ is $\tilde P$-generic over $V$ and 
$G(n)$ is $Q_n$-generic over $V[\tilde G\ast G_n]$.
Fix $\delta<\omega_1$, $n_0\in\omega$ and $f:\omega\to\omega$
and set 
\[
  x=(\eta_n(\alpha))_{n\geq n_0,\delta+f(n)<\alpha<\delta+\omega}.
\]
Then $(G_n)_{n\in\omega}$ is in $V[\tilde G,x]$.

See~Figure~\ref{fig:P}(c).
\begin{proof}
  By induction on $1\leq h\leq n_0$,
  each $\eta_n(\alpha)$ is determined for $n\geq n_0-h$,
  $\delta+\omega\cdot h\leq \alpha<\delta+\omega\cdot (h+1)$.
  By induction on
  limit ordinals $\delta+\omega\cdot n_0<\delta'<\om1$, 
  each $\eta_n(\alpha)$ is determined for 
  $\delta+\omega\cdot n_0<\alpha<\delta'$.
\end{proof}

\begin{remark}\label{rem:note1}
  For all $n_0\in\omega$ and $f:\omega\to\omega$, there are conditions in
  $\tilde P\ast P_\omega$ that determine
  all $\eta_n(m)$ for $n<n_0$ or $m<f(n)$,
  cf.~Figure~\ref{fig:P}(b).
  (The reason is that $\tilde P\ast
  P_{n_0}$ does not add new reals, as we will see in the next lemma, and that
  each
  $Q_n$-condition can be modified at finitely many places.)
  However, these
  conditions are {\em not} dense. (For exactly the same reason: There is a
  condition $p_0$ stating that $(\eta_n(0))_{n\in\omega}$ 
  codes $(\eta_n\restriction \omega)_{n\in\omega}$,
  via a simple injection from $\omega\times\omega$ to
  $\omega$. Then according to
  the last Lemma, no $p'\leq p_0$ can determine all $\eta_n(0)$.)
\end{remark}

\begin{Lem}\label{lem:proper}
  $\tilde P\ast P_n$ forces that 
  $Q_n$ is proper and does not add a new $\omega$-sequence of ordinals.
\end{Lem}

\begin{proof}
  Work in $V'=V[\tilde G\ast G_n]$ and fix some large regular
  cardinal $\chi^*$.

  Let $N^*\esm H^{V'}(\chi^*)$ be a countable elementary submodel
  containing $\tilde G,\eta_{n-1}$ and $p_0\in Q_n$.
  Set $\delta^*=N^*\cap \om1$.
  Let $(D_i)_{i\in\omega}$ list all dense subsets of $Q_n$ that are
  in $N^*$, and assume $D_0=Q_n$. It is enough to show the following:
  \begin{eqpar}\label{eq:XY7}
    There is a $q\leq p_0$ with
    $\height(q)=\delta^*$ such that
    $q$ is stronger than some $p_i\in D_i\cap N^*$
    for every $i\in\omega$.
  \end{eqpar}
  Then $q$ is in particular $N^*$-generic, which shows that
  $Q_n$ is proper. And if $\n f\in N^*$ 
  is a name for a function from
  $\omega$ to the ordinals, then the value of $\n f(n)$ is determined
  in the dense set $D_{i(n)}$ for some $i(n)\in\omega$ 
  an therefore by $q$. This shows that no new $\n f$ is added by $Q_n$.

  So let us prove \ref{eq:XY7}.
  Pick (in $V'$) a sequence $(N_i)_{i\in\omega}$ and a large, regular 
  $\chi\in N^*$ such that:
  \begin{itemize}
    \item $N_i\in N^*$.
    \item $N_i\esm H(\chi)$ is countable.
    \item $N_0$ contains $\eta_{n-1}$ and $p_0$, $N_{i+1}$ contains $N_{i}$
      and $D_{i+1}$.
  \end{itemize}
  Set $\beta_i=N_i\cap\om1$. So $\sup_{i\in\omega}(\beta_i)=\delta^*$.

  Fix (in $V'$) any $\eta^*\in Q_n$ of height $\delta^*$. In particular
  $\eta^*$ coheres with $\eta_{n-1}(\delta^*+m)$ for all $m$.
  Set
  \[
    u_i=\{\alpha\in \nu_{\delta^*,n-1,m}:\, m<i,\beta_i\leq \alpha<\beta_{i+1}\}.
  \]
  Each $u_i$ is finite.
  We will construct $q$ such that $q\supseteq \eta^*\restriction u_i$ for all
  $i\geq 1$.
  This guarantees that $q$ coheres with $\eta_{n-1}(\delta^*+m)$ for all
  $m\in\omega$.

  Assume that $p_i\in N_i\cap D_i$ is already defined.
  We extend it to $p_{i+1}\in N_{i+1}\cap D_{i+1}$:
  The finite sequence $\eta^*\restriction u_{i+1}$ is 
  in $N_{i+1}$, so we can%
  \footnote{by using \ref{lem:basicQ}(i,ii)}
  (in $N_{i+1}$) extend $p_i$ first to some 
  $p'\supseteq \eta^*\restriction u_{i+1}$ in $Q_n$.
  Then extend $p'$ to $p_{i+1}\in D_{i+1}\cap N_{i+1}$.
  
  Set $q=\bigcup_{i\in\omega}p_i$. Then $q$ is in $Q_n$: 
  We already know that $q$ coheres with  $\eta_{n-1}(\delta^*+m)$.
  For $\alpha<\delta^*$,
  let $i$ be such that $\alpha<\beta_i$. Then
  $q$ extends $p_{i+1}$ which coheres with $\eta_{n-1}(\alpha+m)$.
\end{proof}

As an immediate consequence we get the following fact, illustrated
in~Figure~\ref{fig:P}(d):
\begin{Cor}\label{cor:trivial2}
  The conditions $(\tilde p,p_0,\dots,p_{n-1})$ of the following form
  are dense in
  $\tilde P \ast P_n$: $\height(\tilde p)=\delta+\omega\cdot (n-1)$ 
  for some $\delta$, 
  and in $V$ there is a sequence $(p'_0,\dots,p'_{n-1})$ such that
  $p'_i\in 2^{\delta+\omega\cdot (n-1-i)}$ and $p_i$ is the standard name%
  \footnote{
    With ``standard name for $x$'' ($x$ in the ground model)
    we mean the (canonical) name $\check x$ that evaluates to 
    $x$ for all generic filters.
  }
  for $p'_i$.
\end{Cor}

\section{A dense subset}

We will now use the notions of variable, term and substitution as defined in
Definition~\ref{def:terms}. The set of variables we use is $\{x_{n,m}:\,
n,m\in\omega\}$.

Assume that  $\bar p$ is a sequence of terms
$(p_{n,\alpha})_{n\in\omega,\alpha<\delta}$. 
In $V[\tilde G]$, $\bar p$ can be
interpreted as a promise that the generic sequence $(\eta_n)_{n\in \omega}$ is
compatible with $\bar p$, i.e., that 
there is an assignment $a$ such that 
$p_{n,\alpha}\circ a=\eta_n(\alpha)$ for all $n\in\omega,\alpha<\delta$.
Of course such a promise can be inconsistent, for example if 
$\delta=\omega$ and
each $p_{n,m}$ is (the constant term) $0$.

\begin{Def}\label{def:R}
  $R=\bigcup_{\delta<\om1}R_{\delta+\omega}$. 
  A condition $p$ in $R_{\delta+\omega}$ consists of $\tilde p$ and $\bar p$
  such that:
  \begin{itemize}
    \item $\tilde p\in \tilde P$, $\height(\tilde p)=\delta+1$ (or equivalently
      $\delta+\omega$).
    \item $\bar p=(p_{n,\alpha})_{n\in\omega,\alpha<\delta+\omega}$.
    \item $p_{n,\delta+m}$ is the term $x_{n,m}$.
    \item If $\alpha<\delta$, then $p_{n,\alpha}$ is a term that 
      only depends on $x_{l,k}$ with $l<n$.
    \item For every $m,n\in\omega$ and $\alpha\leq \delta$ limit
      there is a $k_0<\omega$ such that
      for all assignments $a$, we get that
      $(p_{n+1,\zeta}\circ a)_{\zeta<\alpha}$ coheres with 
      $p_{n,\alpha+m}\circ a$  above $k_0$.
  \end{itemize}
\end{Def}

We interpret terms are functions, not syntactical objects, so we identify two
elements $p,q$ of $R_{\delta+\omega}$ if they satisfy $\tilde p=\tilde q$ and
$p_{n,\alpha}=^*q_{n,\alpha}$ for all $n,\alpha$; see
Definition~\ref{def:terms}.

Elements of $R$ can be interpreted as statements about the generic
sequence:

\begin{Def}\label{def:acdelta}
  The canonical assignment 
  $a^c_\delta$ assigns the value 
  $\eta_n(\delta+m)$ to the variable $x_{n,m}$.
  (So $a^c_\delta$ is a $\tilde P\ast P_\omega$-name.)
  We also use $a^c_\delta$ as a $\tilde P\ast P_n$-name
  for the partial assignment that maps $\eta_l(\delta+m)$
  to the variable $x_{l,m}$ for all $l<n$.
\end{Def}

\begin{Def}
  Let $i:R\to \tilde P\ast P_\omega$ map $p\in R_{\delta+\omega}$ to
  $(\tilde p,q(0),q(1),\dots)$ defined as follows:
  For each $n$, $q(n)$ is the $\tilde P\ast P_n$-name 
  for the sequence $(p_{n,\alpha}\circ a^c_\delta)_{\alpha<\delta}$.
\end{Def}

\begin{Lem}\label{lem:finitewurscht}
  \begin{enumerate}
   \item $i(p)$ actually is a condition in $\tilde P\ast P_\omega$.
   \item $i(p)$ is the truth value (in $\ro(\tilde P\ast P_\omega$)
     of the following statement:
     $\tilde p\in \tilde G$, and
     $\bar p$ is compatible with the generic sequence $\bar \eta$.
   \item In particular, this truth value is positive.
     Moreover, the truth value remains positive if we additionally assign
     specific values
     for finitely many of the variables $x_{n,m}$.
  \end{enumerate}
\end{Lem}
Here, ``$\bar p$ is compatible with the generic sequence $\bar \eta$'' means:
There is some assignment $a$ such that $p_{n,\alpha}\circ a=\eta_n(\alpha)$ for
all $\alpha< \delta+\omega$. Since $p_{n,\delta+m}=x_{n,m}$, the only
assignment that can ever witness compatibility is the canonical assignment
$a^c_\delta$.

More formally, and slightly stronger, we can formulate the last item as:
Given $f:\omega\to\omega$ and $(b_{n,i})_{n\in \omega,i<f(n)}$ with
$b_{n,i}\in\{0,1\}$, the truth value of the following statement is non-zero:
\begin{itemize}
  \item $\tilde p\in \tilde G$,
  \item $p_{n,\alpha}\circ a^c_\delta=\eta_n(\alpha)$ for all
    $\alpha<\delta+\omega$,
  \item and additionally $\eta_{n}(\delta+i)$ (or equivalently $x_{n,i}\circ
    a^c_\delta$) is  $b_{n,i}$ for all $n\in \omega$ and $i<f(n)$.
\end{itemize}

\begin{proof}
  (i) it follows from the definition  of $R$ that each
  $q(n)$ is a valid condition in $Q_n$.
  (ii) The canonical assignment is the 
  only assignment that can possibly witness compatibility.
  (iii) Given $f$ and $b_{n,i}$ as above, we can just 
  extend $q(n)$ to be the name of some condition $q'$
  in $Q_n$ of height $\delta+\omega$ (instead of just $\delta$)
  such that $q'(\delta+i)=b_{n,i}$ for all $i<f(n)$.
  For this we need, as usual, just Lemma~\ref{lem:basicQ}(i,ii).
\end{proof}

\begin{remark}
  \begin{itemize}
    \item
      It is easy to see (similarly to~\ref{lem:basicQ}) that
      $R_{\delta+\omega}$ is nonempty for all $\delta$. We will
      only prove this (implicitly) for
      ``stationary many''
      $\delta$, in Lemma~\ref{lem:dense}: $\sigma'' R$ is dense in
      $\tilde P\ast P_\omega$.
    \item In view of this Lemma, the proof of
      (iii) can be compared to Remark~\ref{rem:note1}:
      While we cannot densely determine the gray area 
      of~Figure~\ref{fig:P}(b), we can densely determine 
      such an area shifted up to some $\delta$.
  \end{itemize}
\end{remark}

If $\alpha<\delta$, then $p_{n,\alpha}$ can be calculated from finitely many 
$p_{l,\delta+m}$ with $l<n$ (since $p_{n,\alpha}$ is a term using variables
$x_{l,m}$, $l<m$, and  $p_{l,\delta+m}=x_{l,m}$). We can also calculate values
in the other direction:
\begin{Lem}\label{lem:blubb}
  \begin{enumerate}
    \item $p_{n,\alpha}$ is determined by finitely many $x_{l,k}$
     with $l<n$.
    \item
      $x_{n,m}$ 
      can be determined by finitely many $p_{l,k}$ with $l>n,k\in\omega$.
  \end{enumerate}
\end{Lem}

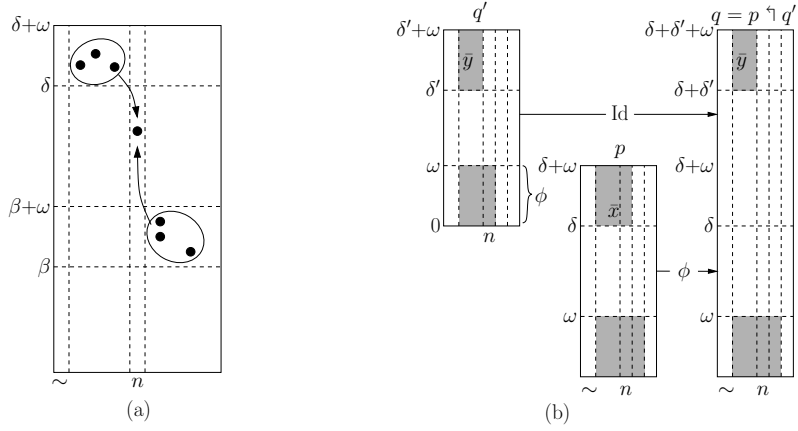
\begin{figure}[tb]
      \newcommand{\mylabelbackandforth}{(a)}
      \newcommand{\mylabelstacking}{(b)}
  \hfill
  \begin{minipage}[t]{.4\textwidth}
    \begin{center}
      \scalebox{0.4}{\input{backandforth.pstex_t}}
    \end{center}
  \end{minipage}
  \hfill
  \begin{minipage}[t]{.55\textwidth}
    \begin{center}
      \scalebox{0.4}{\input{stacking.pstex_t}}
    \end{center}
  \end{minipage}
  \hfill
  \caption{\label{fig:R} Elements of $R$.
     \mylabelbackandforth{}~Dependence in both 
     directions, according to Lemma~\ref{lem:blubb}.
     \mylabelstacking{}~Conditions can be stacked to get stronger
     conditions. If on the other hand $q$ is stronger than
     $p$, then it can be split accordingly.}
\end{figure}

More generally, we get (cf.~Figure~\ref{fig:R}(a)):
If $p\in R_{\delta+\omega}$ and $\beta<\delta$
(not necessary a limit), then
every $p_{n,\alpha}$ with $\beta+\omega\leq \alpha<\delta+\omega$
can be determined by finitely many $p_{l,\zeta}$ with
$l>n,\beta<\zeta<\beta+\omega$. 
More precisely: There is a $k\in\omega$ and a 
sequence $(l_i,\zeta_i)_{i<k}$ 
such that $l_i>n$, $\beta<\zeta_i<\beta+\omega$ and for 
all assignments $a,b$ the following holds:
If $p_{n,\alpha}\circ a\neq p_{n,\alpha}\circ b$, then
$(p_{l_i,\zeta_i}\circ a)_{i< k}\neq (p_{l_i,\zeta_i}\circ b)_{i< k}$.

\begin{proof}
  By induction on $\alpha$:
  Assume $\alpha=\beta+\omega+m$. Then
  $(p_{n+1,\zeta})_{\zeta<\beta+\omega}$ 
  coheres with $p_{n,\alpha}$
  above some $k_0$,
  so we can use $f_{\beta+\omega,n,m}$ to get $p_{n,\alpha}$.
  Now assume that the statement is true for all $\alpha<\nu$, $\nu$ limit.
  If $\alpha=\nu+m$, then $p_{n,\alpha}$ again is determined
  by the values of certain $p_{l,\zeta}$ with $l>n,\beta<\zeta<\nu$,
  each of which in turn is determined (by induction) by finitely many
  $p_{l',\zeta'}$ with $\beta<\zeta'<\beta+\omega$.
\end{proof}

We can identify $R$ with a subset of $P\ast P_\omega$:
\begin{Lem}\label{lem:injective}
  $i:R\to \tilde P\ast P_\omega$ is injective.
\end{Lem}
\begin{proof}
  Fix $p\in R_{\delta+\omega},q\in R_{\delta'+\omega}$, $p\neq q$.
  If $\tilde p\neq \tilde q$, then 
  $i(q)\neq i(p)$. So assume that $\tilde p= \tilde q$ (in particular
  $\delta'=\delta$).
  Since $p\neq q$, there is an $(n,\alpha)$ and a (finite, partial)
  assignment $a$ such that 
  $q_{n,\alpha}\circ a\neq p_{n,\alpha}\circ a$.
  According to \ref{lem:finitewurscht}(iii),
  $i(q)$ is compatible with $a$. Let $r\leq i(q)$ force that 
  the generic sequences are compatible with $a$. Then $r$ forces that 
  $i(p)$ is not in the generic filter, since 
  it determines a different value for $\eta_n(\alpha)$ than $i(q)$.
\end{proof}

So we can interpret $R$ as a subset of $P\ast P_\omega$; and we usually do
so, that is, we will may just write $p$ instead of $i(p)$ and
$R$ instead of $i''R$, as in the following:

\begin{Lem}\label{lem:dense}
  $R\subseteq P\ast P_\omega$ is dense.
\end{Lem}

The proof is a bit cumbersome, but really just a modification of the 
proof of Lemma~\ref{lem:proper}.
\begin{proof}
  Fix $(\tilde p,p(0),p(1),\dots)\in \tilde P\ast P_\omega$,
  and a countable $N^*\esm H(\chi^*)$ containing $(\tilde p,p(0),p(1),\dots)$.
  Set $\delta^*=N^*\cap \om1$. It is enough to show:
  \begin{eqpar}\label{eq:whatweneed}
    There is a $q\in R_{\delta^*+\omega}$ such that
    $i(q)\leq (\tilde p,p(0),p(1),\dots)$.
  \end{eqpar}
  The $\tilde P\ast P_n$-condition $(\tilde p,p(0),\dots,p(n))$
  will be denoted by $p\restriction n$.
  For $q\in R$ and $n\in\omega$, we set
  \[
    q(n)=(q_{n,\alpha})_{\alpha< \delta^*+\omega}\quad\text{and}\quad
    q\restriction n:=(\tilde q,q(0),q(1),\dots,q(n-1))
  \]
  Just as $q$ can be interpreted as
  a condition in $\tilde P\ast P_\omega$ in a canonical way
  (cf.~\ref{lem:finitewurscht}), we can 
  interpret $ q\restriction n$
  as a condition in $\tilde P\ast P_n$. In particular,
  \mbox{``$q\restriction n$ forces $\varphi$''} means the following:
  \begin{quote}
    If $\tilde G\ast G_n$ is $\tilde P\ast P_n$-generic over
    $V$, if $\tilde G$ contains $\tilde q$ and if
    $\eta_0,\dots,\eta_{n-1}$
    are compatible with $(q_{l,\alpha})_{l<n,\alpha< \delta^*+\omega}$,
    then $\varphi$ holds in $V[\tilde G\ast G_n]$.
  \end{quote}

  Let us call an antichain $E$ in $\tilde P\ast P_n$ nice, if 
  every condition $e$ in $E$ has the form of Corollary~\ref{cor:trivial2}. 
  These conditions are dense, so we get:
  \begin{eqpar}\label{eq:nice}
    For all $n\in\omega$, $X\in V$ and all $\n\tau$ such that
    $\tilde P\ast P_n$ forces that $\n\tau\in \check X$ there is a 
    nice maximal antichain $B$ deciding $\n\tau$. I.e., for each
    $b\in B$ there is an $x^b\in X$ such that $b$ forces  $\n\tau=x^b$.
  \end{eqpar}

  {\em The induction hypothesis.}
  We will construct $\tilde q$ in $\tilde Q$ of height $\delta^*+\omega$
  and, by induction on $n\geq 0$,
  the condition $q(n)$ --- i.e., the terms 
  $(q_{n,\alpha})_{\alpha<\delta^*}$ depending on variables
  $x_{i,j}$ with $i<n$ --- such that the following holds:
  \begin{enumerate} 
    \item $q\restriction (n+1)$ satisfies the conditions on elements of
      $R_{\delta^*+\omega}$.
    \item $q\restriction (n+1)$ is $\tilde P\ast P_{n+1}$-generic over $N^*$.
    \item $q\restriction (n+1)$ is stronger than $(\tilde p,p(0),\dots,p(n))$.
    \item $q\restriction (n+1)$ decides every nice maximal antichain $E$ of 
      $\tilde P\ast P_{n+1}$ in $N^*$ by finite case distinction.
  \end{enumerate} 

  More formally: Item (i) means
  \begin{enumerate}
    \item[(i)'] for all $m\in\omega$ and $\alpha\leq \delta^*$ limit
      there is a $k_0<\omega$ such that for all assignments $a$,
      we have that
      $(q_{n,\zeta}\circ a)_{\zeta<\alpha}$ coheres
      with $q_{n-1,\alpha+m}\circ a$
      above $k_0$.
  \end{enumerate} 
  And item (iv) means: For every nice maximal antichain $E$ of 
  $\tilde P\ast P_{n+1}$ in $N^*$ 
  there is an $l^E\in\omega$, a sequence $(e^E_0,\dots e^E_{l^E-1})$ of 
  elements of $E\cap N^*$ and a sequence
  $(t^E_0,\dots,t^E_{l^E-1})$ of terms using 
  only variables $x_{i,j}$ with $i<n+1$
  such that $q\restriction (n+1)$
  forces the following:
  \begin{enumerate}
    \item[(iv)'] There is exactly
      one $k<l^E$ such that 
      $t^E_k\circ a^c_{\delta^*}=1$ (cf.~\ref{def:acdelta}),
      and  $e^E_k\in \tilde G\ast G_n$ for this $k$.
  \end{enumerate}
  This implies the following
  (where we apply Lemma~\ref{lem:finitewurscht}(iii)):
  \begin{enumerate}
     \item[(v)] For all partial assignments $b$ of the
       (finitely many) variables used in any of the $t^E_k$ there is
       exactly one $k<l^E$ such that $t^E_k\circ b=1$.
  \end{enumerate}

  Note that (iii) (for all $n$) implies~\eqref{eq:whatweneed}.

  {\em Step 1: Finding $\tilde q$.}
  First extend $\tilde p$ to 
  $\tilde p'$ such that $\height(\tilde p')=\delta^*$ 
  and such that 
  for every dense subset $D$ of $\tilde P$ in $N^*$
  there is an $d\in D\cap N^*$ weaker that  $\tilde p'$
  (this is possible since $\tilde P$ is $\sigma$-closed).
  In particular, $\tilde p'$ is $\tilde P$-generic over $N^*$,
  and if $E\subseteq \tilde P$ is a maximal antichain
  in $N^*$, then $\tilde p'$ decides the $e\in E$ 
  that will be in the generic filter (and $e\in N^*$).%
  \footnote{
    So  $\tilde p'$ decides ``everything'' about $N^*[\tilde G]$.
    Of course, $N^*[\tilde G]$ is not an element of $V$
    (since it contains, e.g., $\tilde G$). But every formula about 
    $N^*[\tilde G]$ (with parameters in $N^*$)
    is already decided in $V$ ``modulo $\tilde p'$'',
    since every such formula is decided by an antichain.
    We can find such a strong $\tilde p'$ 
    since $\tilde P$ is $\sigma$-complete, and 
    we can do the same for 
    $\tilde P\ast P_1$. However, for $n\geq 1$, $Q_n$
    is not $\sigma$-complete, and we will not be able to decide everything
    with the generic condition $q(n)$;
    but we will still be able to decide ``modulo finite case distinction''.}

  We further extend $\tilde p'$ to $\tilde q$ by adding some arbitrary value
  at $\delta^*$. So $\height(\tilde q)=\delta^*+1$ (or equivalently
  $\delta^*+\omega$).

  {\em Step 2: Finding $q(0)$.}
  This case, $n=1$,
  is simple since $Q_0$ is $\sigma$-closed.

  We have to define the (constant) terms $(q_{0,\alpha})_{\alpha< \delta^*}$.
  Let $(D_i)_{i\in\omega}$ enumerate all $\tilde P$-names
  in $N^*$ for open dense subsets of $Q_0$, such that $D_0=Q_0$

  We now define $r_n$
  and $s_n$ for $n\in\omega$ such that:
  \begin{itemize}
    \item[(a)] $r_n$ is a $\tilde P$-name in $N^*$, forced 
      by $\tilde p$ to be a $Q_0$ condition and element of $D_n$.
    \item[(b)] $s_n$ is a $0$-$1$-sequence in $N^*$, forced by
      $\tilde q$ to be $r_n$.
    \item[(c)] $r_{n+1}$ is forced to extend $s_n$.
  \end{itemize}
  Set $r_0=p(0)$. This satisfies~(a).
  Given an $r_n$ satisfying~(a), note that  
  $\tilde P$ does not add new countable sequences of 
  ordinals. So every condition in $Q_0^{N^*[\tilde G]}$,
  in particular $r_n$,
  already exists in the ground model $N^*$. So $r_n$ 
  is decided by a  maximal antichain, and therefore 
  by $\tilde q$, to be some sequence $s_n\in N^*$; satisfying~(b).
  Also, since $s_n\in 2^{<\om1}\cap N^*$,  we can find in $N^*$
  a $\tilde P$-name $r_{n+1}$ 
  for an element of $D_{n+1}$ extending $s_n$.

  Fix $\alpha<\delta^*$, and set
  $q_{0,\alpha}$ to be the term with constant value $s_n(\alpha)$
  (for sufficiently large $n$). This defines $q(0)$.
  So $\tilde q$ forces that that $q(0)$ is $Q_0$-generic over
  $N^*[\tilde G]$,
  i.e., $(\tilde q, q(0))$ is $\tilde P\ast P_1$-generic 
  over $N^*$ and forces 
  that $(\tilde p,p(0))\in \tilde G\ast G_1$. So 
  (i)--(iii) are satisfied.
  Now fix some nice, maximal antichain 
  $E\subset \tilde P\ast P_1$ such that $E\in N^*$. 
  Every $e\in E$ is of the form $(\tilde e, e(0))$
  for a $0$-$1$-sequence $e(0)$ in $V$.
  If $e\in N^*$, then $\tilde e$  and $e(0)$ have height
  less than $\delta^*$.
  In particular, every $e\in E\cap N^*$ is either extended by 
  $(\tilde q,q(0))$  or is incompatible with it.
  Since $E\in N^*$ is a maximal antichain, and
  since $(\tilde q,q(0))$ is $\tilde P\ast P_1$-generic over 
  $N^*$, we know that there has to be 
  exactly one $e^*\in E$ compatible with $(\tilde q,q(0))$,
  and $e^*\in N^*$.
  In other words, $(\tilde q,q(0))$ decides the
  element $e^*\in E\cap N^*$ that is going to be in 
  $\tilde G\ast G_1$. So to satisfy (iv)', we can set
  $l^E=1$, $t^E_0=1$, $e^E_0=e^*$. 

  {\em Step 3: The successor step.}
  Now things get a bit more complicated, since $Q_{n}$ 
  is not $\sigma$-closed any more.
  We assume that the induction hypothesis (i)--(iv)
  holds for $n-1$. So we already have 
  $q\restriction n$
  want to find $q_{n,\alpha}$ for $0\leq \alpha<\delta^*$.
  As previously, we let $(D_i)_{i\in\omega}$ enumerate all 
  $\tilde P\ast P_n$-names in $N^*$ for open dense subsets of $Q_n$
  (and we set $D_0=Q_n$).

  First we fix (in $V$) a term-sequence
  $(t^*_\alpha)_{\alpha\in \bigcup_{m\in\omega}\nu_{\delta^*,n-1,m}}$ 
  such that:
  \begin{itemize}
    \item If $\alpha\in \nu_{\delta^*,n-1,m}$, then $t^*_\alpha$ only depends
      on $x_{n-1,m}$.
    \item For all $m\in\omega$, the sequence $\bar t^*$ coheres with
      $q_{n-1,\delta^*+m}$ (which is just $x_{n-1,m}$)
      above some $k_0$.
    \item For every $\beta<\delta^*$, the partial sequence
      $\bar t^*\restriction \beta$ uses only finitely many variables.
  \end{itemize}
  We can find such a sequence since the $f_{\delta^*,n-1,m,k}$ defined
  by $\tilde q$ are surjective and the $\nu_{\delta^*,n-1,m}$ are disjoint
  (for different $m$) cofinal subsets of $\delta^*$ of order type $\omega$.

  We will construct in $V$ by induction on $i\in\omega$
  \begin{itemize}
    \item a finite set $v_i$
      of variables $x_{l,j}$ with $l<n$,
    \item for every (partial) assignment $a$ of $v_i$ a
      $\tilde P\ast P_n$-name $r^a_i$ in $N^*$,
    \item a finite set $w_i$
      of variables $x_{l,j}$ with $l<n$, 
    \item 
      for every assignment $b$ of $w_i$ a
      $0$-$1$-sequence $s^b_i$ in $N^*$,
    \item an ordinal $\beta_i<\delta^*$,
  \end{itemize}
  such that the following holds:
  \begin{itemize}
    \item[(a)] $v_{i+1}\supseteq w_i\supseteq v_i$.
    \item[(b)] If $a$ is an assignment of $v_i$, then
      $r^a_i$ is a name (in $N^*$) for an element of $D_i$.
    \item[(c)] If $b$ is an assignment of $w_i$ and
      $a$ its restriction to $v_i$, then
      $(q\restriction n)\&b$ forces%
      \footnote{
        For $x\in \tilde P\ast P_n$,
        $x\& b$ is the 
        truth value (in $\ro(\tilde P\ast P_n)$)
        of the following statement:
        $x\in \tilde G*G_n$, and $b$ is compatible
        with $\eta_0,\dots,\eta_{n-1}$, i.e.,
        for every $x_{l,k}\in w_i$, we have 
        $\eta_{l,\delta^*+k}=x_{l,k}\circ b$. For this notation
        we can 
	use $x=p\restriction n$, and also $x=q\restriction n$,
        since we can canonically interpret
        $q\restriction n$ as element of $\tilde P\ast P_n$.%
      }
      $s^b_i=r^a_i$.
  \end{itemize}

  Set $w_{0}=\emptyset$. So there is only one assignment,
  the empty one, of $w_0$. We set $r^\emptyset_0=p(n)$.%
        \footnote{More formally, we should set
          \[
            r^\emptyset_0=
            \begin{cases}
              p(n)&\text{if }p(n)\in Q_n\\
              \emptyset   &\text{otherwise,}
            \end{cases}
          \]
          since $p(n)$ is forced to be in $Q_n$ by $p\restriction n$,
          not by the empty condition.}
  Assume that for some $i\geq 0$ we already have 
  $w_i$, and $r^a_i$ for all assignments $a$ of $w_i$.  Fix $a$.
  Note that $\tilde P\ast P_n$ does not add any new countable
  sequences of ordinals, so according to~\eqref{eq:nice}
  $r^a_i$ is decided by a nice maximal antichain $E$
  of $\tilde P\ast P_n$ in $N^*$.
  Using item (iv) of the induction hypothesis, we choose
  the sequences $(e^E_0,\dots e^E_{l^E})$ of
  and $(t^E_0,\dots,t^E_{l^E})$.
  Let $v'$ be the (finite) set of variables used
  in any of the $t^E_k$. Set $v^a_{i}=w_{i}\cup v'$.
  Let $b$ be an assignment of $v^a_i$ extending $a$.
  According to (v), there is a unique $k\leq l^E$ such that
  $t^{E}_k\circ b=1$. We call this element $k(b)$.
  The element $e^E_{k(b)}$ determines $r^a_i$ to be a 
  specific $0$-$1$-sequence of $V$, and we call this sequence $s^b_i$.
  Note that $s^b_i\in N^*$. We can do this for all assignments
  $a$ of $w_i$, and set $v_i=\bigcup v^a_i$.

  We still have to construct $\beta_i$, $w_{i+1}$ and $r^b_{i+1}$.
  We pick in $N^*$ a $\tilde P\ast P_n$-name $N$ for a 
  countable elementary submodel of
  $H^{V[\tilde G \ast G_n]}(\chi)$
  containing $D_{i+1}$ and all the (finitely many) $s^b_i$. Since 
  $N\cap\om1$ is an $\tilde P\ast P_n$-name for an
  ordinal, 
  there are only finitely many possibilities modulo $q\restriction n$,
  and we can choose $\beta_{i}\in N^*\cap \om1$ 
  larger than every possibility for $N\cap\om1$.
  
  The terms $t^*_\alpha$ for $\alpha<\delta_i$ use only 
  a finite set $w'$ of
  variables (of the form $x_{n-1,m}$). Set $w_{i+1}=v_i\cup w'$.
  Fix an assignment $a$ of $w_{i+1}$ and let $b$ be the restriction
  to $v_i$. Fix the index set
  \[
    I=\{\alpha<\beta_i:\, (\exists j\leq i) \alpha\in \nu_{\delta^*,n-1,j}\}.
  \]
  The finite set $I$ is in $N^*$. Set
  \[
    \bar x=(\bar t^*\circ b)\restriction I.
  \]
  This is a finite partial function in $N^*$ from $I$ to $\{0,1\}$.
  We define the $\tilde P\ast P_n$-name
  $r^a_{i+1}$ in $N^*$ by the following construction
  in $N^*[\tilde G\ast G_n]$: 
  (Let $d$ be some fixed element of $D_{i+1}$.)
  \begin{itemize}
    \item Assume that 
      $\beta'=N\cap \omega_1<\beta_{i}$. (Otherwise set $r^a_{i+1}=d$.)
    \item Assume that $s^b_i$ is a $Q_n$-condition.
      (Otherwise set $r^a_{i+1}=d$.)
    \item In $N$, extend $s^b_i$ to some $Q_n$-condition 
      containing
      $\bar x\restriction (\beta'\setminus \height(s^b_i))$.
      (As usual, use~\ref{lem:basicQ} inside $N$.)
    \item Again in $N$, pick some condition $r^a_{i+1}$ in
      $D_{i+1}$ extending $s'$.  In particular,
      $r_{i+1}$ has height less than $\beta_i$.
  \end{itemize}

  This ends the construction.
  We can summarize all the possibilities of 
  $s^a_i(\alpha)$ into the term $q_{n,\alpha}$ (depending on 
  the variables in $v_i$).
  This defines $q(n)$. 

  It remains to be shown that $q\restriction n+1$ satisfies the
  induction hypothesis.

  For~(i)', first assume $\alpha<\delta^*$.
  Let $D_i$ be the set of conditions of length $\geq\alpha$.
  Let $q_{n-1,\alpha+m}$ be determined by the finite set $v$
  of variables, and set $v'=v\cup v_i$.
  Fix an assignment $b$ of $v'$. In particular $b$ determines 
  $q_{n-1,\alpha+m}$ as well as $q(n)\restriction \alpha$,
  since $q(n)$ ``extends'' $s^{b'}_i\in D_i$
  (where $b'$ is the restriction of $b$ to $v_i$).
  Since $q\restriction n$ is compatible with the finite assignment
  $b$, we know that  $q(n)\circ b\restriction \alpha$ coheres
  with $q_{n-1,\alpha+m}\circ b$ above some $k^b_0$.
  So we can set $k_0$ to be the maximum of all the $k^b_0$
  for all assignments $b$ of $v'$.

  Now assume $\alpha=\delta^*$ and $m\in \omega$.
  Pick $\gamma\in \nu_{\delta^*,n-1,m}\setminus \beta_{m}$.
  Look at the term $q_{n,\gamma}$. According to the 
  construction, 
  \[
    q_{n,\gamma}\circ b=s^{b'}_i(\gamma)\circ b=t^*_\alpha\circ b
  \]
  for all assignments $b$, and therefore $q(n)$ coheres with
  $x_{n-1,m}$.

  Let us now show~(iv).
  Let $E\in N^*$ be a nice, maximal antichain of $\tilde P\ast P_{n+1}$.
  Let $D$ be the $\tilde P\ast P_{n}$-name for the following open
  dense subset of $Q_n$
  \[
    D=\{q\leq e(n):\, e\in E,e\restriction n\in \tilde G\ast G_n\}.
  \]
  We know that $D$ appears as some $D_i$ in the list of dense sets in $N^*$.
  Fixing an assignment $a$ of $v_i$, we get
  $s^a_i$ in $N^*$ such that $s^a_i\in D_i[\tilde G\ast G_n]$. We set
  \[
    A^a=\{e\restriction n:\, e\in E, e(n)\subseteq s^a_i \}
  \]
  This is a nice $\tilde P\ast P_{n}$-antichain and maximal under
  $(q\restriction n)\& a$. We can extend it to a nice maximal antichain
  $B^a$. By induction hypothesis, we can determine 
  modulo $q\restriction n$ the element $b$ of $B^a$ chosen by 
  $\tilde G\ast G_n$ filter by finite case distinction.
  Then $b^\frown s^a_i$ is the element of $E$ chosen by $\tilde G\ast G_{n+1}$.
  Combining the finite case distinction for the $s^a_i$ with the
  finite case distinctions for the according $B^a$ gives the desired result.
\end{proof}

Since $R$ is a subset of $\tilde P \ast P_\omega$, it is also a partial order
(and since it is dense, it is forcing equivalent to $\tilde P \ast P_\omega$).
We now show that we can interpret the order on $R$ in a different way,
using substitutions of terms:

\begin{Def}
  Let $p\in R_{\delta+\omega}$ and $q\in R_{\delta'+\omega}$. We call $q$
  term-stronger than $p$, if either $p=q$ or if the following holds:
  $\tilde q\leq \tilde p$ (in particular $\delta'\geq \delta$), and
  $p_{n,\alpha}\circ \phi =^* q_{n,\alpha}$ for all $\alpha<\delta$ and
  for the substitution $\phi$ defined by $\phi_{n,m}=q_{n,\delta+m}$.
\end{Def}
(Again, recall that we interpret terms as functions, so we use $=^*$ as defined
in~\ref{def:terms}.)

\begin{Lem}
  The condition $q\in R$ is term-stronger than $p\in R$ iff $i(q)\leq i(p)$.
\end{Lem}

\begin{proof}
  Assume that $q$ is not term-stronger than $p$.
  If $\tilde q$ is not stronger than $\tilde p$ in $\tilde P$,
  then $i(q)$ cannot be stronger than $i(p)$. So assume $\tilde q\leq \tilde p$.
  According to the definition of term-stronger,
  $q_{l,\alpha}=^* p_{l,\alpha}\circ \phi$ fails for some $l,\alpha$.
  These terms depend on finitely many variables $x_{n,m}$, and 
  there is a partial assignment $a$ of these variables
  such that $q_{l,\alpha}\circ a\neq p_{l,\alpha}\circ \phi\circ a$.
  According to
  Lemma~\ref{lem:finitewurscht}(iii), 
  we can force the generic sequence to be compatible with 
  $q$ and $a$. Then $i(q)$
  is in the generic filter, but $i(p)$ is not, contradicting
  $i(q)\leq i(q)$.
\end{proof}


%

If $q'$ is a condition, then $(q'_{n,m})_{n,m\in\omega}$ can be 
interpreted as substitution: 
For $p\in R_{\delta+\omega}$
and  $q'\in R_{\delta'+\omega}$,
we can stack $q'$ on top of $p$
--- overlapping at $[\delta,\delta+\omega[$ --- to get a condition
$q \in R_{\delta+\delta'+\omega}$ stronger than $p$,
cf.~Figure~\ref{fig:R}(b). We write $q=p\Lsh q'$.

More precisely:
\begin{Def}\label{def:stack}
  For $p\in R_{\delta+\omega}$ and $q'\in R_{\delta'+\omega}$,
  we define the condition $q=p\Lsh q'$ in $R_{\delta+\delta'+\omega}$ 
  as follows 
  \begin{itemize}
    \item $\tilde q\restriction(\delta+1)=\tilde p$.
    \item $\tilde q(\alpha)$ for $\alpha\geq \delta+\omega$ is defined 
      the following way:
      \\
        $\nu^q_{\delta+\alpha,n,m}=
        \{\delta+\beta:\, \beta\in \nu^{q'}_{\alpha,n,m}\}$,
      \\
        $j^q_{\delta+\alpha,n,m}=j^{q'}_{\alpha,n,m}$,
      \\
        $f^q_{\delta+\alpha,n,m,k}=f^{q'}_{\alpha,n,m,k}$.
    \item $q_{n,\delta+\alpha}=q'_{n,\alpha}$.
    \item If $\alpha<\delta$, then $q_{n,\alpha}=p_{n,\alpha}\circ \phi$ for 
      the substitution $\phi$ defined by $\phi_{n,m}=q'_{n,m}$.
  \end{itemize}
\end{Def}

\begin{figure}[tb]
  \newcommand{\labelstrangea}{(a)}
  \newcommand{\labelstrangeb}{(b)}
  \hfill
  \begin{minipage}[t]{.4\textwidth}
    \begin{center}
      \scalebox{0.4}{\input{strangea.pstex_t}}
    \end{center}
  \end{minipage}
  \hfill
  \begin{minipage}[t]{.4\textwidth}
    \begin{center}
      \scalebox{0.4}{\input{strangeb.pstex_t}}
    \end{center}
  \end{minipage}
  \hfill
  \caption{\label{fig:strange} \labelstrangea{} If $q$ is stronger than
    $p1$ and $p_2$, then $q_{2,m}$ is constant for all $m$. (The
    gray area indicates constant terms; the $\tilde P$ parts are not
    displayed.) \labelstrangeb{} If $p_1,p_2$ are compatible (i.e., 
    weaker than some $q$), 
    they do not have to be comparable.}
\end{figure}
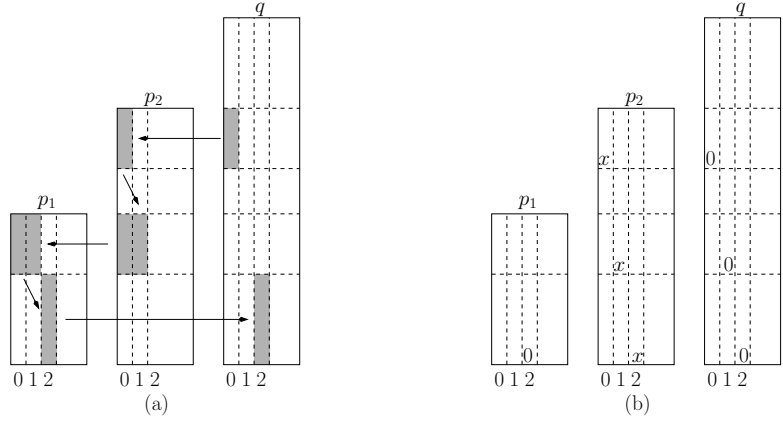

\begin{Fact}
  \begin{enumerate}
    \item
      If $p\in R_{\delta+\omega}$ and $q'\in R_{\delta'+\omega}$,
      then $p\Lsh q'$ is stronger than $p$.
    \item
      If $q\in R_{\delta+\delta'+\omega}$ is stronger than $p\in
      R_{\delta+\omega}$, then we can ``split'' $q$ into $p\in
      R_{\delta+\omega}$ and $q'\in R_{\delta'+\omega}$ such that $q=p\Lsh q'$.
  \end{enumerate}
\end{Fact}
\begin{remarks}
  \begin{itemize}
    \item
      Of course we generally cannot split a condition at every level: If $ q\in
      R_{\delta''+\omega}$ and $\delta'<\delta''$, then we generally do not
      get $q=p\Lsh q'$ for some $p\in R_{\delta'+\omega}$.
    \item The Fact shows that for all $q\in R$ there are only finitely many
      $p\geq q$, see Figure~\ref{fig:strange}(a).
    \item Note that two compatible conditions generally are not comparable,
      see Figure~\ref{fig:strange}(b). (Otherwise, according to
      the previous item, $R$ would be isomorphic to
      a tree of height $\omega$ and therefore collapse the continuum.)
    \item 
      The situation is similar to $Q_*$ defined in
      Section~\ref{sec:sacks1}: The conditions that are stronger than some $p\in
      R$ are exactly those with another condition $q'\in R$ stacked on top.
  \end{itemize}
\end{remarks}

\section{Sacks reals as squares of terms again}\label{sec:notagain}

We will now investigate the relation of $R$ and $Q_*$.  Given a $p\in R$, we
can restrict $p$ to an $\omega\times \omega$-matrix of terms:
\begin{Def}
  For $p\in R$, set $\sigma(p)=(p_{n,m})_{n,m\in\omega}$.
\end{Def}

Note that
\begin{equation}
  \sigma(p\Lsh q')=\sigma(p)\circ \sigma(q').
\end{equation}
So stacking $q'$ on top of $p$ translates to applying $\sigma(q')$ (as
substitution) to $\sigma(p)$.

Generally $\sigma(p)$ will not be element of $Q_*$, and for a $\bar t\in
Q_*$ there generally is no $p\in R$ such that $\sigma(p)=\bar t$. The
reason is that some obvious conditions on the term-matrix are incomparable: In
$Q_*$, we require
\begin{quote}
  $t_{i,j}$ only depends on $x_{n,m}$ such that $(n,m)\unlhd(i,j)$,
\end{quote}
whereas every $\sigma(p)$ obviously satisfies
\begin{quote}
  $p_{i,j}$ only depends on $x_{n,m}$ such that $n<i$.
\end{quote}

We will now define a dense subset $R'\subseteq R$ such that
$\sigma''R'\subseteq Q_*$, and such that $R'$ adds a $Q_*$-generic object. This
proves the first part of Theorem~\ref{thm:main}, since $Q_*$ is forcing
equivalent to Sacks forcing and $R'$ is (as a dense subset) equivalent to $R$,
which in turn is dense in $P*P_\omega$. So $P*P_\omega$ adds a Sacks real.

\begin{Lem}\label{lem:laestig}
  \begin{itemize}
    \item[(i)] There is an $\rdelta\in R_{\omega\cdot\omega+\omega}$
      such that $\sigma(\rdelta)\in Q_*$
      and  $\sigma(p\Lsh \rdelta)\in Q_*$ for all $p\in R$.
    \item[(ii)] There is an $\rmulti\in R_{\omega\cdot \omega+\omega}$ 
      such that for all $\phi\in Q_*$ there is an $r\in R_{\omega\cdot
      \omega+\omega}$ with 
      $\sigma(\rmulti)\circ  \phi=\sigma(r)$.
  \end{itemize}
\end{Lem}

We postpone the proof to the end of the section.
We set 
      \[
        R'=\{p\Lsh \rdelta:\, p\in R\}.
      \]
This is a dense subset of $R$, since $p\Lsh \rdelta\leq p$ for all $p$.
As a consequence of the previous Lemma, we get:
\begin{Cor}\label{cor:gurke2}
  \begin{enumerate}
    \item[(a)] If $l\in R'$ and $p\leq_R l$, then
      $\sigma(p\Lsh \rdelta)\leq_{Q_*} \sigma(l)$.
    \item[(b)] If $p\in R'$ then there is an  $\bar s\leq\sigma(p)$
      such that for all $\bar t\leq \bar s$  then there is 
      an $q\leq p$ in $R'$ such that $\sigma(q)\leq \bar t$.
    \item[(c)] The forcing notion $R'$ adds a generic for $Q_*$.
      So $\tilde P\ast P_\omega$ adds a Sacks real.
  \end{enumerate}
\end{Cor}

\begin{proof}[Proof of the Corollary]
  (a)
  Assume that $p=l\Lsh q'$. Then $\sigma(p\Lsh \rdelta)=\sigma(l)\circ
  \sigma(q'\Lsh \rdelta)$; and $\sigma(q'\Lsh \rdelta)$ is an element of $Q_*$
  and therefore witnesses that $\sigma(p\Lsh \rdelta)$ is
  stronger than $\sigma(l)$.

  (b) Set
       $ \bar s=\sigma(p\Lsh \rmulti\Lsh \rdelta)$, and
      let $\phi$ witnesses $\bar t\leq \bar s$, i.e., $\phi\in Q_*$ and
      \[
	\bar t=\bar s\circ \phi=\sigma(p)\circ \sigma(\rmulti)\circ
        \sigma(\rdelta)\circ \phi.
      \]
      $ \sigma(\rdelta)\circ \phi\in Q_*$, so by
      Lemma~\ref{lem:laestig}(ii),
      there is an $r\in R$ such that 
      \[
        \sigma(r)=\sigma(\rmulti)\circ \sigma(\rdelta)\circ \phi,\text{ so }
      \]
      \[
        \sigma(p\Lsh r)=\sigma(p)\circ \sigma(r)=\bar s\circ \phi=\bar t.
      \]
      Set $q=p \Lsh r\Lsh \rdelta$. Then  $q\in R'$ and $q\leq p$.
      Furthermore, 
      \[
        \sigma (q)=\sigma (p\Lsh r)\circ \sigma(\rdelta)=\bar t\circ
        \sigma(\rdelta)\leq \bar t.
      \]
  (c) Let $G'$ be $R'$-generic over $V$. We show
  that the 
  following set is $Q_*$-generic filter over $V$:
  \begin{equation}\label{eq:brubra}
    G_*=\{\bar r\in Q_*:\, (\exists q \in G')\, \sigma(q)\leq \bar r\}
  \end{equation}
  First note that $ G_{*}$ does not
  contain incompatible elements: Assume that $\bar r_1$ and $\bar r_2$
  are in $G_*$. Then there are $l_1,l_2\in G'$ such that $\sigma(l_i)\leq r_i$.
  Since $G'$ is a filter, there is some $p\leq l_1,l_2$ in $G'$.
  The set 
  \[
    \{p'\Lsh \rdelta:\, p'\leq p\}
  \]
  is dense below $p$, so there is some $q=p'\Lsh \rdelta$ in $G'$.
  According to (a), the $q$
  satisfies $\sigma(q)\leq \sigma(l_1),\sigma(l_2)$.

  Now assume that $D\subseteq Q_{*}$ is dense, 
  and (towards a contradiction) that $p$ forces that $G_{*}$ does not 
  meet $D$. Then pick $\bar s$ as in (b), pick 
  $\bar t\leq \bar s$ in $D$ and pick
  $q$ again as in (b). So
  $q$ forces that $\sigma(p)\leq t$ is in $G_*$, a
  contradiction. So we know that~\eqref{eq:brubra} is generic.
\end{proof}

It remains to prove Lemma~\ref{lem:laestig}.
All these facts are easy to {\em see}, but a bit cumbersome to write down 
formally. So the reader might be better off drawing a picture than
reading the proof.
Fix an injective function from $\omega^{<\omega}$ to $\omega$,
$(a_1,\dots,a_l)\mapsto \ulcorner a_1,\dots,a_l\urcorner$, with coinfinite
range. 

\subsection*{The construction of $\rdelta$}
All we need is a $\rdelta\in R_{\omega\cdot\omega+\omega}$ satisfying
the following:
\begin{eqpar}\label{eq:hert23}
  $\rdelta_{n,m}$ only depends on variables $x_{i,j}$ such that
  $i+j<n$.
\end{eqpar}
Then, if we stack $\rdelta$ on top of any $p\in R_{\delta+\omega}$, 
the resulting $q=p\Lsh
\rdelta$
will satisfy~\eqref{eq:hert23} as well.
Also, every element $q$ of $R$ satisfies that each $x_{i,j}$ depends on
finitely many $q_{n,m}$ for $n,m\in\omega$, according to
Lemma~\ref{lem:blubb}(ii).  Therefore $\sigma(p\Lsh \rdelta)$ will
satisfy all requirements for an element of $Q_*$, which proves
Lemma~\ref{lem:laestig}(i).

We now construct $\rdelta$.
\begin{itemize}
  \item When defining $\tilde\rdelta$,
    only the $\nu$ part is nontrivial; we set each
    $j_{\alpha,n,m}:\omega\to \omega$ and $f_{\alpha,n,m,k}:2\to 2$ 
    to be the
    identity function for all $\alpha,n,m,k$.%
    \footnote{This corresponds to the
      simpler version of $\tilde P$ in Section~\ref{sec:easy}.%
    }
  \item 
    We deal with one variable at a time. Assume that we deal
    with $x_{n_0,m_0}$.
  \item Set $\nu_{\omega\cdot \omega,n_0,m_0}=
    \{\omega\cdot k+\ulcorner n_0,m_0\urcorner:\, k>m_0\}$.\\ 
    For $\alpha\in \nu_{\omega\cdot \omega,n_0,m_0}$, we set
    $\rdelta_{n_0+1,\alpha}=x_{n_0,m_0}$.
  \item
    If $n=n_0+l$ for some $l\geq 1$, if  $m=\ulcorner
    a_0,a_1,\dots,a_l\urcorner$  with $a_0=n_0,a_1=m_0$, and if $k>m_0-l$,
    then set
    $\nu_{\omega\cdot (k+1),n, m }= 
    \{\omega\cdot k+\ulcorner a_0,\dots,a_l,j\urcorner :\, j\in\omega\}$,
    and
    for $\alpha\in \nu_{\omega\cdot (k+1),n, m }$, we set
    $\rdelta_{n+1,\alpha}=x_{n_0,m_0}$.
  \item We repeat this for all $x_{i,j}$. (Note that the 
    $\nu_{\alpha,n, m }$ defined for different $m$ will be disjoint).
  \item So far, whenever we have defined some $\nu_{\beta,n,m}$
    to contain $\alpha$, we also guaranteed that 
    $\rdelta_{n+1,\alpha}$ and $\rdelta_{n,\beta+m}$ are the same
    variable. 
  \item
    We now set all $\rdelta_{n,\alpha}$ that are undefined so far
    to be the constant term $0$, and define every
    $\nu_{\beta,n,m}$ that is undefined so far in a way such that 
    every member $\alpha$ of $\nu_{\beta,n,m}$ 
    satisfies $\rdelta_{n+1,\alpha}=0$. (Here, we use that the coding
    function has coinfinite range.)
\end{itemize}
It is easy to see that the object $\rdelta$ defines this way is element of $R$.
Each $\rdelta_{n,\alpha}$ is either an $x_{n_0,m_0}$ or $0$.  If
$\rdelta_{n,m}=x_{n_0,m_0}$, then $n>m_0+n_0$. Given $n$, $\rdelta_{n,m}=0$ for
infinitely many $m$.

\subsection*{The construction of $\rmulti$}
We will first show the following:
\begin{Lem}
  If $(p_{n,m})_{n,m\in\omega}$ satisfies
      \begin{itemize}
        \item[1.] $p_{n,m}$ is a term depending only on $x_{i,j}$ with $i<n$,
        \item[2.] $(\forall n)\,(\exists^\infty m)\,p_{n,m}=0$,
	\item[3.] $(\forall i,j)\,(\exists^\infty n)\,(\forall
              M)\,(\exists m_0,\dots m_k>M)\, 
              x_{i,j}$ is determined by $p_{n,m_0},\dots,p_{n,m_k}$,
      \end{itemize}
      then there is a $q\in R_{\omega\cdot \omega+\omega}$ such that 
      $\sigma(q)=\bar p$.
\end{Lem}

\begin{proof}
The proof is very similar to the preceding construction.
The reader might just consult Figure~\ref{fig:simple}(c).

Assume we have such a sequence $\bar p$.  We have to define $q\in
R_{\omega\cdot \omega+\omega}$. We already know that
$q_{n,m}=p_{n,m}$ for $n,m\in\omega$.

We more or less repeat the construction above, to get all $q_{n,\alpha}$ and
all $\nu_{\beta,n,m}$, but only for $\beta\geq \omega+\omega$,
and we deal
with $\nu_{\omega,n,m}$ later.
Assume we are dealing with $x_{n_0,m_0}$.
Set 
\[
  I_{n_0,m_0}=\{n'\geq n_0+2:\, (\forall M)\,(\exists m'_0,\dots
  m'_k>M)\,x_{n_0,m_0} \text{ is determined by }p_{n',m'_0},\dots,p_{n',m'_k}\}.
\]
According to assumption~(3), $I_{n_0,m_0}$ is infinite.

For all $\alpha>\omega$ limit and all $n,m,k$ we will set $j_{\alpha,n,m}$ and
$f_{\alpha,n,m,k}$ to be the identity functions.

We set 
\[
  \nu_{\omega\cdot\omega,n_0,m_0}=\{i\cdot\omega + \ulcorner n_0,m_0 \urcorner:\, i+m_0+1\in I_{n_0,m_0} \}.
\]
(So $\nu_{\omega\cdot\omega,n_0,m_0}\cap [i\cdot \omega,(i+1)\cdot \omega[$
contains a singleton if $i+m_0+1\in I_{n_0,m_0}$, and is empty otherwise.)
We set $q_{n_0+1,\alpha}=x_{n_0,m_0}$ for all $\alpha\in \nu_{\omega\cdot\omega,n_0,m_0}$;
and ``propagate the $x_{n_0,m_0}$ diagonally down'': If $n=n_0+l$
for some $l\geq 1$, if $i>0$ and $i+m+l\in I_{n_0,m_0}$, and if $m=\ulcorner
a_0,a_1,\dots,a_l\urcorner$ such that $a_0=n_0,a_1=m_0$,
then set 
\[
  \nu_{(i+1)\cdot \omega,n,m}=\{i\cdot\omega + \ulcorner a_0,\dots,a_l,j
  \urcorner: j\in\omega\}.
\]

We iterate this for all $x_{n_0,m_0}$, and set
all $q_{n,\alpha}$ that have not been defined in this process to be
the constant term $0$.
Also we set the $\nu_{\beta,n,m}$ for $\beta>\omega$
that have not been defined yet to
contain only $\alpha>\omega$ such that $q_{n+1,\alpha}=0$.
(Remember that the coding had coinfinite range.)

So we have all $q_{n,\alpha}$ and all $\nu_{\beta,n,m}$,
$j_{\beta,n,m}$ and $f_{\beta,n,m,k}$ for $\beta>\omega$.

We still have to define $\nu_{\omega,n,m}$, $j_{\omega,n,m}$ and 
$f_{\omega,n,m,k}$. For this, we use a simple book-keeping:
At stage $i$, there are only finitely many pairs $(n,m)$
for which any of these objects are already partially defined.
For all of these $(n,m)$, we also have:
$\nu_{\omega,n,m}$ is defined up to height $M_{n,m}$,
$j_{\omega,n,m}$ is defined up to some $h_{n,m}$ such that
$j_{\omega,n,m}(h_{n,m}-1)=M_{n,m}-1$  
$f_{\omega,n,m,k}$ is defined for exactly the $k<h_{n,m}$.
Let $M$ be the maximum of all $M_{n,m}$ for a given stage.

The book-keeping gives us an $(n_0,m_0)$ and an $(n,m)$ such that
$q_{n,\omega+m}=x_{n_0,m_0}$. By our construction, we know that $x_{n_0,m_0}$
can be determined by finitely many and arbitrary large $q_{n+1,m'}$. Fix
$m'_0,\dots, m'_{l-1}$ bigger than $M$ such that $q_{n+1,m'_0},\dots
,q_{n+1,m'_{l-1}}$ determines $x_{n_0,m_0}$.  Extend $\nu_{\omega,n,m'}$ to
contain exactly $\{m'_0,\dots, m'_l\}$, continue $j_{\omega,n,m'}$ by setting
$j_{\omega,n,m'}(h_{n,m})=M_{n,m}+l-1$ and define $f_{\omega,n,m',h_{n,m}}$  so
that it calculates $x_{n_0,m_0}$.

At the end, again set the $\nu_{\omega,n,m}$ that have not been defined
in this process to contain only $m'$ such that $q_{n+1,m'}=0$.
To be able to do this, we use at height $\omega$ assumption~(2).
\end{proof}

We can now define $\rmulti$: We can take any condition in $R$ satisfying
\begin{itemize}
  \item $\rmulti_{n,m}$ only depends on $x_{i,j}$ with $i+j<n$.
  \item Every $x_{i,j}$ with $i+j<n$, as well as the constant $0$ term,
    occurs infinitely often in $\{\rmulti_{n,m}:\, m\in\omega\}$.
\end{itemize}
If we set $\bar p= \rmulti\circ \phi$ for some $\phi\in Q_*$, we get:
\begin{itemize}
  \item $p_{n,m}$ only depends on $x_{i,j}$ with $i+j<n$. (Due to 
    \eqref{eq:blug}.) So we satisfy (1).
  \item For all $n$, infinitely many $p_{n,m}$ are $0$. So we satisfy (2).
  \item $x_{i,j}$ is determined by $(\phi_{l_i,k_i})_{i\in I}$.
    Fix any $n$ bigger than $\max(l_i+k_i:\, i\in I)$. Then
    $x_{i,j}$ is determined by finitely many $p_{n,m}$ (where
    we can pick the $m$'s arbitrarily large). So we satisfy (3).
\end{itemize}
So $\rmulti\circ \phi$ satisfies all assumptions of the previous Lemma, and we
get a $q$ as desired.

\section{The quotient forcing.}\label{sec:quotient}

It might look tempting to assume $\Diamond$ to construct the coding
sequences $\bar \nu,\bar j,\bar f$ instead of using the preparatory forcing. 
(We just have to ``guess'' correctly sufficiently often for the
proofs to work.) However, this is not possible: 
Otherwise, Sacks forcing would be {\em equivalent} to $P_\omega$ (since
$P_\omega$ adds a Sacks real $s$ which in turn determines the
$P_\omega$-generic filter $G_\omega$). But Sacks reals are minimal, and the
$Q_0$ generic $\eta_0\in 2^{\om1}$ is not in the ground model $V$.  Therefore
$V[s]=V[\eta_0]$, a contradiction to the fact that $V[\eta_0]$ does not
add new reals.

In particular,
if we look at $P_\omega$ in $V[\tilde G]$,
then $P_\omega$ does not add  a Sacks real (over $V[\tilde G]$), 
just a Sacks real over $V$.

So $\tilde P\ast P_\omega$ adds a Sacks real $s$ but is not equivalent to Sacks
forcing, and $s$ does not determine the $\tilde P$-generic object  $\tilde G$.
However, every new $\omega$-sequence is already added by $s$:
\begin{Lem}
  If $\tilde G\ast G_\omega$ is $\tilde P\ast P_\omega$-generic over
  $V$, and if $r\in V[\tilde G\ast G_\omega]$ is an $\omega$-sequence of
  ordinals, then $r\in V[s]$. Here we set
  $s=(\bar \eta_n(m))_{n,m\in\omega}$, the Sacks real over $V$.
\end{Lem}
\begin{proof}
  If $\tilde q\in \tilde G$ has height $\delta$, then
  $s$ together with $\tilde q$ determines $G_n$ up to height $\delta$
  for all $n$
  (just as in Lemma~\ref{lem:newreal}).
  So if $N\esm H(\chi)$ and $\tilde q\in \tilde G$ has height $N\cap \om1$, then
  $s$ together with $\tilde q$ determines 
  whether $r\in \tilde G\ast G_\omega$ for any $r\in R\cap N$.

  Assume towards a contradiction 
  that $p\in R$ forces that $\n f$ is an $\omega$-sequence
  of ordinals not added by $s$. Choose an $N\esm H(\chi)$ containing 
  $p,\n f$, and an $N$-generic $q\leq p$. 
  Each $\n f(n)$ is decided by some maximal antichain $A\in N$.
  But for each $a\in A\cap N$,
  $s$ together with $\tilde q$ determines whether $a$ is in $G$.
  In particular, $\n f[G]\in V[s]$.
\end{proof}

This proves the second part of Theorem~\ref{thm:main}:
Since $R$ forces that there is some Sacks real over $V$ and
since Sacks forcing is homogeneous, $R$ can be factored 
as Sacks composed with some $P'$. Since the Sacks real
already adds all new $\omega$-sequences, $P'$ is NNR.

\bibliographystyle{plain}
\bibliography{905}
\end{document}

%% file: tree.pstex_t
\begin{picture}(0,0)%
\includegraphics{tree.pstex}%
\end{picture}%
\setlength{\unitlength}{4144sp}%
\begingroup\makeatletter\ifx\SetFigFontNFSS\undefined%
\gdef\SetFigFontNFSS#1#2#3#4#5{%
  \reset@font\fontsize{#1}{#2pt}%
  \fontfamily{#3}\fontseries{#4}\fontshape{#5}%
  \selectfont}%
\fi\endgroup%
\begin{picture}(3942,1572)(-3884,-6121)
\put(-3756,-5194){\makebox(0,0)[lb]{\smash{{\SetFigFontNFSS{14}{16.8}{\rmdefault}{\mddefault}{\updefault}{\color[rgb]{0,0,0}$F_2$}%
}}}}
\put(-3869,-4786){\makebox(0,0)[b]{\smash{{\SetFigFontNFSS{14}{16.8}{\rmdefault}{\mddefault}{\updefault}{\color[rgb]{0,0,0}$T$}%
}}}}
\put(-674,-4786){\makebox(0,0)[b]{\smash{{\SetFigFontNFSS{14}{16.8}{\rmdefault}{\mddefault}{\updefault}{\color[rgb]{0,0,0}$T'$}%
}}}}
\put(-224,-5236){\makebox(0,0)[lb]{\smash{{\SetFigFontNFSS{14}{16.8}{\rmdefault}{\mddefault}{\updefault}{\color[rgb]{0,0,0}$F'_0$}%
}}}}
\put(-3644,-5700){\makebox(0,0)[lb]{\smash{{\SetFigFontNFSS{14}{16.8}{\rmdefault}{\mddefault}{\updefault}{\color[rgb]{0,0,0}$F_1$}%
}}}}
\put(-2384,-6001){\makebox(0,0)[lb]{\smash{{\SetFigFontNFSS{14}{16.8}{\rmdefault}{\mddefault}{\updefault}{\color[rgb]{0,0,0}$F_0$}%
}}}}
\end{picture}%

%% file: order.pstex_t
\begin{picture}(0,0)%
\includegraphics{order.pstex}%
\end{picture}%
\setlength{\unitlength}{4144sp}%
\begingroup\makeatletter\ifx\SetFigFontNFSS\undefined%
\gdef\SetFigFontNFSS#1#2#3#4#5{%
  \reset@font\fontsize{#1}{#2pt}%
  \fontfamily{#3}\fontseries{#4}\fontshape{#5}%
  \selectfont}%
\fi\endgroup%
\begin{picture}(2558,1381)(2067,-4861)
\put(3418,-3886){\makebox(0,0)[rb]{\smash{{\SetFigFontNFSS{12}{14.4}{\rmdefault}{\mddefault}{\updefault}{\color[rgb]{0,0,0}$3$}%
}}}}
\put(3418,-4111){\makebox(0,0)[rb]{\smash{{\SetFigFontNFSS{12}{14.4}{\rmdefault}{\mddefault}{\updefault}{\color[rgb]{0,0,0}$2$}%
}}}}
\put(3418,-4336){\makebox(0,0)[rb]{\smash{{\SetFigFontNFSS{12}{14.4}{\rmdefault}{\mddefault}{\updefault}{\color[rgb]{0,0,0}$1$}%
}}}}
\put(3418,-4561){\makebox(0,0)[rb]{\smash{{\SetFigFontNFSS{12}{14.4}{\rmdefault}{\mddefault}{\updefault}{\color[rgb]{0,0,0}$0$}%
}}}}
\put(3587,-4786){\makebox(0,0)[b]{\smash{{\SetFigFontNFSS{12}{14.4}{\rmdefault}{\mddefault}{\updefault}{\color[rgb]{0,0,0}$0$}%
}}}}
\put(3812,-4786){\makebox(0,0)[b]{\smash{{\SetFigFontNFSS{12}{14.4}{\rmdefault}{\mddefault}{\updefault}{\color[rgb]{0,0,0}$1$}%
}}}}
\put(4037,-4786){\makebox(0,0)[b]{\smash{{\SetFigFontNFSS{12}{14.4}{\rmdefault}{\mddefault}{\updefault}{\color[rgb]{0,0,0}$2$}%
}}}}
\put(4262,-4786){\makebox(0,0)[b]{\smash{{\SetFigFontNFSS{12}{14.4}{\rmdefault}{\mddefault}{\updefault}{\color[rgb]{0,0,0}$3$}%
}}}}
\put(2926,-4786){\makebox(0,0)[b]{\smash{{\SetFigFontNFSS{12}{14.4}{\rmdefault}{\mddefault}{\updefault}{\color[rgb]{0,0,0}$3$}%
}}}}
\put(2701,-4786){\makebox(0,0)[b]{\smash{{\SetFigFontNFSS{12}{14.4}{\rmdefault}{\mddefault}{\updefault}{\color[rgb]{0,0,0}$2$}%
}}}}
\put(2476,-4786){\makebox(0,0)[b]{\smash{{\SetFigFontNFSS{12}{14.4}{\rmdefault}{\mddefault}{\updefault}{\color[rgb]{0,0,0}$1$}%
}}}}
\put(2251,-4786){\makebox(0,0)[b]{\smash{{\SetFigFontNFSS{12}{14.4}{\rmdefault}{\mddefault}{\updefault}{\color[rgb]{0,0,0}$0$}%
}}}}
\put(2082,-4561){\makebox(0,0)[rb]{\smash{{\SetFigFontNFSS{12}{14.4}{\rmdefault}{\mddefault}{\updefault}{\color[rgb]{0,0,0}$0$}%
}}}}
\put(2082,-4336){\makebox(0,0)[rb]{\smash{{\SetFigFontNFSS{12}{14.4}{\rmdefault}{\mddefault}{\updefault}{\color[rgb]{0,0,0}$1$}%
}}}}
\put(2082,-4111){\makebox(0,0)[rb]{\smash{{\SetFigFontNFSS{12}{14.4}{\rmdefault}{\mddefault}{\updefault}{\color[rgb]{0,0,0}$2$}%
}}}}
\put(2082,-3886){\makebox(0,0)[rb]{\smash{{\SetFigFontNFSS{12}{14.4}{\rmdefault}{\mddefault}{\updefault}{\color[rgb]{0,0,0}$3$}%
}}}}
\put(2476,-4111){\makebox(0,0)[b]{\smash{{\SetFigFontNFSS{12}{14.4}{\rmdefault}{\mddefault}{\updefault}{\color[rgb]{0,0,0}$7$}%
}}}}
\put(3590,-4561){\makebox(0,0)[b]{\smash{{\SetFigFontNFSS{12}{14.4}{\rmdefault}{\mddefault}{\updefault}{\color[rgb]{0,0,0}$0$}%
}}}}
\put(3590,-4561){\makebox(0,0)[b]{\smash{{\SetFigFontNFSS{12}{14.4}{\rmdefault}{\mddefault}{\updefault}{\color[rgb]{0,0,0}$0$}%
}}}}
\put(3815,-4561){\makebox(0,0)[b]{\smash{{\SetFigFontNFSS{12}{14.4}{\rmdefault}{\mddefault}{\updefault}{\color[rgb]{0,0,0}$2$}%
}}}}
\put(3590,-4336){\makebox(0,0)[b]{\smash{{\SetFigFontNFSS{12}{14.4}{\rmdefault}{\mddefault}{\updefault}{\color[rgb]{0,0,0}$1$}%
}}}}
\put(3590,-4111){\makebox(0,0)[b]{\smash{{\SetFigFontNFSS{12}{14.4}{\rmdefault}{\mddefault}{\updefault}{\color[rgb]{0,0,0}$3$}%
}}}}
\put(3590,-3886){\makebox(0,0)[b]{\smash{{\SetFigFontNFSS{12}{14.4}{\rmdefault}{\mddefault}{\updefault}{\color[rgb]{0,0,0}$6$}%
}}}}
\put(3815,-4336){\makebox(0,0)[b]{\smash{{\SetFigFontNFSS{12}{14.4}{\rmdefault}{\mddefault}{\updefault}{\color[rgb]{0,0,0}$4$}%
}}}}
\put(3815,-4111){\makebox(0,0)[b]{\smash{{\SetFigFontNFSS{12}{14.4}{\rmdefault}{\mddefault}{\updefault}{\color[rgb]{0,0,0}$7$}%
}}}}
\put(4040,-4561){\makebox(0,0)[b]{\smash{{\SetFigFontNFSS{12}{14.4}{\rmdefault}{\mddefault}{\updefault}{\color[rgb]{0,0,0}$5$}%
}}}}
\put(4040,-4336){\makebox(0,0)[b]{\smash{{\SetFigFontNFSS{12}{14.4}{\rmdefault}{\mddefault}{\updefault}{\color[rgb]{0,0,0}$8$}%
}}}}
\put(4265,-4561){\makebox(0,0)[b]{\smash{{\SetFigFontNFSS{12}{14.4}{\rmdefault}{\mddefault}{\updefault}{\color[rgb]{0,0,0}$9$}%
}}}}
\put(2251,-4561){\makebox(0,0)[b]{\smash{{\SetFigFontNFSS{12}{14.4}{\rmdefault}{\mddefault}{\updefault}{\color[rgb]{0,0,0}$0$}%
}}}}
\put(2251,-4561){\makebox(0,0)[b]{\smash{{\SetFigFontNFSS{12}{14.4}{\rmdefault}{\mddefault}{\updefault}{\color[rgb]{0,0,0}$0$}%
}}}}
\put(2476,-4561){\makebox(0,0)[b]{\smash{{\SetFigFontNFSS{12}{14.4}{\rmdefault}{\mddefault}{\updefault}{\color[rgb]{0,0,0}$2$}%
}}}}
\put(2251,-4336){\makebox(0,0)[b]{\smash{{\SetFigFontNFSS{12}{14.4}{\rmdefault}{\mddefault}{\updefault}{\color[rgb]{0,0,0}$1$}%
}}}}
\put(2251,-4111){\makebox(0,0)[b]{\smash{{\SetFigFontNFSS{12}{14.4}{\rmdefault}{\mddefault}{\updefault}{\color[rgb]{0,0,0}$3$}%
}}}}
\put(2251,-3886){\makebox(0,0)[b]{\smash{{\SetFigFontNFSS{12}{14.4}{\rmdefault}{\mddefault}{\updefault}{\color[rgb]{0,0,0}$6$}%
}}}}
\put(2476,-4336){\makebox(0,0)[b]{\smash{{\SetFigFontNFSS{12}{14.4}{\rmdefault}{\mddefault}{\updefault}{\color[rgb]{0,0,0}$4$}%
}}}}
\put(2701,-4561){\makebox(0,0)[b]{\smash{{\SetFigFontNFSS{12}{14.4}{\rmdefault}{\mddefault}{\updefault}{\color[rgb]{0,0,0}$5$}%
}}}}
\put(2701,-4336){\makebox(0,0)[b]{\smash{{\SetFigFontNFSS{12}{14.4}{\rmdefault}{\mddefault}{\updefault}{\color[rgb]{0,0,0}$8$}%
}}}}
\put(2926,-4561){\makebox(0,0)[b]{\smash{{\SetFigFontNFSS{12}{14.4}{\rmdefault}{\mddefault}{\updefault}{\color[rgb]{0,0,0}$9$}%
}}}}
\end{picture}%

%% file: simplecoherence.pstex_t
\begin{picture}(0,0)%
\includegraphics{simplecoherence.pstex}%
\end{picture}%
\setlength{\unitlength}{4144sp}%
\begingroup\makeatletter\ifx\SetFigFontNFSS\undefined%
\gdef\SetFigFontNFSS#1#2#3#4#5{%
  \reset@font\fontsize{#1}{#2pt}%
  \fontfamily{#3}\fontseries{#4}\fontshape{#5}%
  \selectfont}%
\fi\endgroup%
\begin{picture}(2547,5942)(841,-8016)
\put(856,-3301){\makebox(0,0)[rb]{\smash{{\SetFigFontNFSS{20}{24.0}{\rmdefault}{\mddefault}{\updefault}{\color[rgb]{0,0,0}$\delta$}%
}}}}
\put(856,-2761){\makebox(0,0)[rb]{\smash{{\SetFigFontNFSS{20}{24.0}{\rmdefault}{\mddefault}{\updefault}{\color[rgb]{0,0,0}$\delta+m$}%
}}}}
\put(2034,-4426){\makebox(0,0)[b]{\smash{{\SetFigFontNFSS{20}{24.0}{\rmdefault}{\mddefault}{\updefault}{\color[rgb]{0,0,0}$0$}%
}}}}
\put(2034,-5506){\makebox(0,0)[b]{\smash{{\SetFigFontNFSS{20}{24.0}{\rmdefault}{\mddefault}{\updefault}{\color[rgb]{0,0,0}$0$}%
}}}}
\put(2034,-5776){\makebox(0,0)[b]{\smash{{\SetFigFontNFSS{20}{24.0}{\rmdefault}{\mddefault}{\updefault}{\color[rgb]{0,0,0}$1$}%
}}}}
\put(2034,-6316){\makebox(0,0)[b]{\smash{{\SetFigFontNFSS{20}{24.0}{\rmdefault}{\mddefault}{\updefault}{\color[rgb]{0,0,0}$0$}%
}}}}
\put(2031,-4156){\makebox(0,0)[b]{\smash{{\SetFigFontNFSS{20}{24.0}{\rmdefault}{\mddefault}{\updefault}{\color[rgb]{0,0,0}$0$}%
}}}}
\put(1584,-2761){\makebox(0,0)[b]{\smash{{\SetFigFontNFSS{20}{24.0}{\rmdefault}{\mddefault}{\updefault}{\color[rgb]{0,0,0}$0$}%
}}}}
\put(2034,-6046){\makebox(0,0)[b]{\smash{{\SetFigFontNFSS{20}{24.0}{\rmdefault}{\mddefault}{\updefault}{\color[rgb]{0,0,0}$0$}%
}}}}
\put(1576,-7531){\makebox(0,0)[b]{\smash{{\SetFigFontNFSS{20}{24.0}{\rmdefault}{\mddefault}{\updefault}{\color[rgb]{0,0,0}$\eta_{n\mathord-1}$}%
}}}}
\put(2156,-7911){\makebox(0,0)[b]{\smash{{\SetFigFontNFSS{20}{24.0}{\rmdefault}{\mddefault}{\updefault}{\color[rgb]{0,0,0}\mylabelsimplecoherence}%
}}}}
\put(2034,-4966){\makebox(0,0)[b]{\smash{{\SetFigFontNFSS{20}{24.0}{\rmdefault}{\mddefault}{\updefault}{\color[rgb]{0,0,0}$0$}%
}}}}
\put(2036,-3886){\makebox(0,0)[b]{\smash{{\SetFigFontNFSS{20}{24.0}{\rmdefault}{\mddefault}{\updefault}{\color[rgb]{0,0,0}$0$}%
}}}}
\put(2034,-4696){\makebox(0,0)[b]{\smash{{\SetFigFontNFSS{20}{24.0}{\rmdefault}{\mddefault}{\updefault}{\color[rgb]{0,0,0}$1$}%
}}}}
\put(2034,-5236){\makebox(0,0)[b]{\smash{{\SetFigFontNFSS{20}{24.0}{\rmdefault}{\mddefault}{\updefault}{\color[rgb]{0,0,0}$0$}%
}}}}
\put(2026,-7531){\makebox(0,0)[b]{\smash{{\SetFigFontNFSS{20}{24.0}{\rmdefault}{\mddefault}{\updefault}{\color[rgb]{0,0,0}$q$}%
}}}}
\end{picture}%

%% file: R.pstex_t
\begin{picture}(0,0)%
\includegraphics{R.pstex}%
\end{picture}%
\setlength{\unitlength}{4144sp}%
\begingroup\makeatletter\ifx\SetFigFontNFSS\undefined%
\gdef\SetFigFontNFSS#1#2#3#4#5{%
  \reset@font\fontsize{#1}{#2pt}%
  \fontfamily{#3}\fontseries{#4}\fontshape{#5}%
  \selectfont}%
\fi\endgroup%
\begin{picture}(3222,6062)(841,-8016)
\put(1136,-7521){\makebox(0,0)[b]{\smash{{\SetFigFontNFSS{20}{24.0}{\rmdefault}{\mddefault}{\updefault}{\color[rgb]{0,0,0}$\sim$}%
}}}}
\put(856,-2221){\makebox(0,0)[rb]{\smash{{\SetFigFontNFSS{20}{24.0}{\rmdefault}{\mddefault}{\updefault}{\color[rgb]{0,0,0}$\delta\mathord+\omega$}%
}}}}
\put(856,-4246){\makebox(0,0)[rb]{\smash{{\SetFigFontNFSS{20}{24.0}{\rmdefault}{\mddefault}{\updefault}{\color[rgb]{0,0,0}$\delta$}%
}}}}
\put(1396,-4021){\makebox(0,0)[lb]{\smash{{\SetFigFontNFSS{20}{24.0}{\rmdefault}{\mddefault}{\updefault}{\color[rgb]{0,0,0}$x_{0,0}$}%
}}}}
\put(1396,-3706){\makebox(0,0)[lb]{\smash{{\SetFigFontNFSS{20}{24.0}{\rmdefault}{\mddefault}{\updefault}{\color[rgb]{0,0,0}$x_{0,1}$}%
}}}}
\put(3376,-7521){\makebox(0,0)[b]{\smash{{\SetFigFontNFSS{20}{24.0}{\rmdefault}{\mddefault}{\updefault}{\color[rgb]{0,0,0}$n$}%
}}}}
\put(2476,-7911){\makebox(0,0)[b]{\smash{{\SetFigFontNFSS{20}{24.0}{\rmdefault}{\mddefault}{\updefault}{\color[rgb]{0,0,0}\mylabelR}%
}}}}
\put(3106,-3706){\makebox(0,0)[rb]{\smash{{\SetFigFontNFSS{20}{24.0}{\rmdefault}{\mddefault}{\updefault}{\color[rgb]{0,0,0}$x_{n\mathord-1,1}$}%
}}}}
\put(3106,-4021){\makebox(0,0)[rb]{\smash{{\SetFigFontNFSS{20}{24.0}{\rmdefault}{\mddefault}{\updefault}{\color[rgb]{0,0,0}$x_{n\mathord-1,0}$}%
}}}}
\put(1891,-3706){\makebox(0,0)[lb]{\smash{{\SetFigFontNFSS{20}{24.0}{\rmdefault}{\mddefault}{\updefault}{\color[rgb]{0,0,0}$\cdots$}%
}}}}
\put(1891,-4021){\makebox(0,0)[lb]{\smash{{\SetFigFontNFSS{20}{24.0}{\rmdefault}{\mddefault}{\updefault}{\color[rgb]{0,0,0}$\cdots$}%
}}}}
\put(1576,-3436){\makebox(0,0)[b]{\smash{{\SetFigFontNFSS{20}{24.0}{\rmdefault}{\mddefault}{\updefault}{\color[rgb]{0,0,0}$\vdots$}%
}}}}
\put(2926,-3436){\makebox(0,0)[b]{\smash{{\SetFigFontNFSS{20}{24.0}{\rmdefault}{\mddefault}{\updefault}{\color[rgb]{0,0,0}$\vdots$}%
}}}}
\put(856,-6316){\makebox(0,0)[rb]{\smash{{\SetFigFontNFSS{20}{24.0}{\rmdefault}{\mddefault}{\updefault}{\color[rgb]{0,0,0}$\alpha$}%
}}}}
\put(3376,-6136){\makebox(0,0)[b]{\smash{{\SetFigFontNFSS{20}{24.0}{\rmdefault}{\mddefault}{\updefault}{\color[rgb]{0,0,0}$p_{n,\alpha}$}%
}}}}
\end{picture}%

%% file: simpleexample.pstex_t
\begin{picture}(0,0)%
\includegraphics{simpleexample.pstex}%
\end{picture}%
\setlength{\unitlength}{4144sp}%
\begingroup\makeatletter\ifx\SetFigFontNFSS\undefined%
\gdef\SetFigFontNFSS#1#2#3#4#5{%
  \reset@font\fontsize{#1}{#2pt}%
  \fontfamily{#3}\fontseries{#4}\fontshape{#5}%
  \selectfont}%
\fi\endgroup%
\begin{picture}(4122,6107)(166,-8016)
\put(2268,-7911){\makebox(0,0)[b]{\smash{{\SetFigFontNFSS{20}{24.0}{\rmdefault}{\mddefault}{\updefault}{\color[rgb]{0,0,0}\mylabelsimpleexample}%
}}}}
\put(451,-7531){\makebox(0,0)[b]{\smash{{\SetFigFontNFSS{20}{24.0}{\rmdefault}{\mddefault}{\updefault}{\color[rgb]{0,0,0}$\sim$}%
}}}}
\put(181,-6451){\makebox(0,0)[rb]{\smash{{\SetFigFontNFSS{20}{24.0}{\rmdefault}{\mddefault}{\updefault}{\color[rgb]{0,0,0}$\omega$}%
}}}}
\put(181,-3121){\makebox(0,0)[rb]{\smash{{\SetFigFontNFSS{20}{24.0}{\rmdefault}{\mddefault}{\updefault}{\color[rgb]{0,0,0}$\omega\mathord\cdot\omega$}%
}}}}
\put(181,-4651){\makebox(0,0)[rb]{\smash{{\SetFigFontNFSS{20}{24.0}{\rmdefault}{\mddefault}{\updefault}{\color[rgb]{0,0,0}$\omega\mathord\cdot 3$}%
}}}}
\put(181,-5551){\makebox(0,0)[rb]{\smash{{\SetFigFontNFSS{20}{24.0}{\rmdefault}{\mddefault}{\updefault}{\color[rgb]{0,0,0}$\omega\mathord\cdot 2$}%
}}}}
\put(181,-2176){\makebox(0,0)[rb]{\smash{{\SetFigFontNFSS{20}{24.0}{\rmdefault}{\mddefault}{\updefault}{\color[rgb]{0,0,0}$\omega\mathord\cdot\omega\mathord+\omega$}%
}}}}
\put(181,-2671){\makebox(0,0)[rb]{\smash{{\SetFigFontNFSS{20}{24.0}{\rmdefault}{\mddefault}{\updefault}{\color[rgb]{0,0,0}$\omega\mathord\cdot\omega\mathord+2$}%
}}}}
\put(3376,-7531){\makebox(0,0)[b]{\smash{{\SetFigFontNFSS{20}{24.0}{\rmdefault}{\mddefault}{\updefault}{\color[rgb]{0,0,0}$8$}%
}}}}
\put(2926,-7531){\makebox(0,0)[b]{\smash{{\SetFigFontNFSS{20}{24.0}{\rmdefault}{\mddefault}{\updefault}{\color[rgb]{0,0,0}$7$}%
}}}}
\put(2476,-7531){\makebox(0,0)[b]{\smash{{\SetFigFontNFSS{20}{24.0}{\rmdefault}{\mddefault}{\updefault}{\color[rgb]{0,0,0}$6$}%
}}}}
\put(2026,-7531){\makebox(0,0)[b]{\smash{{\SetFigFontNFSS{20}{24.0}{\rmdefault}{\mddefault}{\updefault}{\color[rgb]{0,0,0}$5$}%
}}}}
\put(1576,-7531){\makebox(0,0)[b]{\smash{{\SetFigFontNFSS{20}{24.0}{\rmdefault}{\mddefault}{\updefault}{\color[rgb]{0,0,0}$4$}%
}}}}
\end{picture}%

%% file: coherence.pstex_t
\begin{picture}(0,0)%
\includegraphics{coherence.pstex}%
\end{picture}%
\setlength{\unitlength}{4144sp}%
\begingroup\makeatletter\ifx\SetFigFontNFSS\undefined%
\gdef\SetFigFontNFSS#1#2#3#4#5{%
  \reset@font\fontsize{#1}{#2pt}%
  \fontfamily{#3}\fontseries{#4}\fontshape{#5}%
  \selectfont}%
\fi\endgroup%
\begin{picture}(2547,5942)(841,-8016)
\put(2473,-4259){\makebox(0,0)[lb]{\smash{{\SetFigFontNFSS{20}{24.0}{\rmdefault}{\mddefault}{\updefault}{\color[rgb]{0,0,0}$0$}%
}}}}
\put(2481,-3886){\makebox(0,0)[lb]{\smash{{\SetFigFontNFSS{20}{24.0}{\rmdefault}{\mddefault}{\updefault}{\color[rgb]{0,0,0}$0$}%
}}}}
\put(2480,-6174){\makebox(0,0)[lb]{\smash{{\SetFigFontNFSS{20}{24.0}{\rmdefault}{\mddefault}{\updefault}{\color[rgb]{0,0,0}$1=g_0$}%
}}}}
\put(2473,-5776){\makebox(0,0)[lb]{\smash{{\SetFigFontNFSS{20}{24.0}{\rmdefault}{\mddefault}{\updefault}{\color[rgb]{0,0,0}$1=g_1$}%
}}}}
\put(2475,-5236){\makebox(0,0)[lb]{\smash{{\SetFigFontNFSS{20}{24.0}{\rmdefault}{\mddefault}{\updefault}{\color[rgb]{0,0,0}$0=g_2$}%
}}}}
\put(2386,-6856){\makebox(0,0)[b]{\smash{{\SetFigFontNFSS{20}{24.0}{\rmdefault}{\mddefault}{\updefault}{\color[rgb]{0,0,0}$f$}%
}}}}
\put(2034,-6046){\makebox(0,0)[b]{\smash{{\SetFigFontNFSS{20}{24.0}{\rmdefault}{\mddefault}{\updefault}{\color[rgb]{0,0,0}$0$}%
}}}}
\put(2026,-7525){\makebox(0,0)[b]{\smash{{\SetFigFontNFSS{20}{24.0}{\rmdefault}{\mddefault}{\updefault}{\color[rgb]{0,0,0}$n$}%
}}}}
\put(1584,-2761){\makebox(0,0)[b]{\smash{{\SetFigFontNFSS{20}{24.0}{\rmdefault}{\mddefault}{\updefault}{\color[rgb]{0,0,0}$0$}%
}}}}
\put(2031,-4156){\makebox(0,0)[b]{\smash{{\SetFigFontNFSS{20}{24.0}{\rmdefault}{\mddefault}{\updefault}{\color[rgb]{0,0,0}$0$}%
}}}}
\put(2036,-3886){\makebox(0,0)[b]{\smash{{\SetFigFontNFSS{20}{24.0}{\rmdefault}{\mddefault}{\updefault}{\color[rgb]{0,0,0}$1$}%
}}}}
\put(2034,-4966){\makebox(0,0)[b]{\smash{{\SetFigFontNFSS{20}{24.0}{\rmdefault}{\mddefault}{\updefault}{\color[rgb]{0,0,0}$1$}%
}}}}
\put(2034,-6316){\makebox(0,0)[b]{\smash{{\SetFigFontNFSS{20}{24.0}{\rmdefault}{\mddefault}{\updefault}{\color[rgb]{0,0,0}$0$}%
}}}}
\put(2034,-5776){\makebox(0,0)[b]{\smash{{\SetFigFontNFSS{20}{24.0}{\rmdefault}{\mddefault}{\updefault}{\color[rgb]{0,0,0}$1$}%
}}}}
\put(2034,-5506){\makebox(0,0)[b]{\smash{{\SetFigFontNFSS{20}{24.0}{\rmdefault}{\mddefault}{\updefault}{\color[rgb]{0,0,0}$0$}%
}}}}
\put(2034,-4426){\makebox(0,0)[b]{\smash{{\SetFigFontNFSS{20}{24.0}{\rmdefault}{\mddefault}{\updefault}{\color[rgb]{0,0,0}$0$}%
}}}}
\put(2156,-7911){\makebox(0,0)[b]{\smash{{\SetFigFontNFSS{20}{24.0}{\rmdefault}{\mddefault}{\updefault}{\color[rgb]{0,0,0}\mylabelcoherence}%
}}}}
\put(856,-2761){\makebox(0,0)[rb]{\smash{{\SetFigFontNFSS{20}{24.0}{\rmdefault}{\mddefault}{\updefault}{\color[rgb]{0,0,0}$\delta+m$}%
}}}}
\put(856,-3301){\makebox(0,0)[rb]{\smash{{\SetFigFontNFSS{20}{24.0}{\rmdefault}{\mddefault}{\updefault}{\color[rgb]{0,0,0}$\delta$}%
}}}}
\put(1576,-7531){\makebox(0,0)[b]{\smash{{\SetFigFontNFSS{20}{24.0}{\rmdefault}{\mddefault}{\updefault}{\color[rgb]{0,0,0}$n\mathord-1$}%
}}}}
\end{picture}%

%% file: possible.pstex_t
\begin{picture}(0,0)%
\includegraphics{possible.pstex}%
\end{picture}%
\setlength{\unitlength}{4144sp}%
\begingroup\makeatletter\ifx\SetFigFontNFSS\undefined%
\gdef\SetFigFontNFSS#1#2#3#4#5{%
  \reset@font\fontsize{#1}{#2pt}%
  \fontfamily{#3}\fontseries{#4}\fontshape{#5}%
  \selectfont}%
\fi\endgroup%
\begin{picture}(2547,5942)(841,-8016)
\put(2156,-7911){\makebox(0,0)[b]{\smash{{\SetFigFontNFSS{20}{24.0}{\rmdefault}{\mddefault}{\updefault}{\color[rgb]{0,0,0}\mylabelpossible}%
}}}}
\put(856,-4876){\makebox(0,0)[rb]{\smash{{\SetFigFontNFSS{20}{24.0}{\rmdefault}{\mddefault}{\updefault}{\color[rgb]{0,0,0}$\omega$}%
}}}}
\put(1666,-7486){\makebox(0,0)[b]{\smash{{\SetFigFontNFSS{20}{24.0}{\rmdefault}{\mddefault}{\updefault}{\color[rgb]{0,0,0}$n_0$}%
}}}}
\put(991,-7486){\makebox(0,0)[b]{\smash{{\SetFigFontNFSS{20}{24.0}{\rmdefault}{\mddefault}{\updefault}{\color[rgb]{0,0,0}$\sim$}%
}}}}
\end{picture}%

%% file: impossible.pstex_t
\begin{picture}(0,0)%
\includegraphics{impossible.pstex}%
\end{picture}%
\setlength{\unitlength}{4144sp}%
\begingroup\makeatletter\ifx\SetFigFontNFSS\undefined%
\gdef\SetFigFontNFSS#1#2#3#4#5{%
  \reset@font\fontsize{#1}{#2pt}%
  \fontfamily{#3}\fontseries{#4}\fontshape{#5}%
  \selectfont}%
\fi\endgroup%
\begin{picture}(2547,5942)(841,-8016)
\put(2156,-7911){\makebox(0,0)[b]{\smash{{\SetFigFontNFSS{20}{24.0}{\rmdefault}{\mddefault}{\updefault}{\color[rgb]{0,0,0}\mylabelimpossible}%
}}}}
\put(1666,-7486){\makebox(0,0)[b]{\smash{{\SetFigFontNFSS{20}{24.0}{\rmdefault}{\mddefault}{\updefault}{\color[rgb]{0,0,0}$n_0$}%
}}}}
\put(991,-7486){\makebox(0,0)[b]{\smash{{\SetFigFontNFSS{20}{24.0}{\rmdefault}{\mddefault}{\updefault}{\color[rgb]{0,0,0}$\sim$}%
}}}}
\put(856,-3571){\makebox(0,0)[rb]{\smash{{\SetFigFontNFSS{20}{24.0}{\rmdefault}{\mddefault}{\updefault}{\color[rgb]{0,0,0}$\delta\mathord+\omega$}%
}}}}
\put(856,-5776){\makebox(0,0)[rb]{\smash{{\SetFigFontNFSS{20}{24.0}{\rmdefault}{\mddefault}{\updefault}{\color[rgb]{0,0,0}$\delta$}%
}}}}
\end{picture}%

%% file: typicalfinite.pstex_t
\begin{picture}(0,0)%
\includegraphics{typicalfinite.pstex}%
\end{picture}%
\setlength{\unitlength}{4144sp}%
\begingroup\makeatletter\ifx\SetFigFontNFSS\undefined%
\gdef\SetFigFontNFSS#1#2#3#4#5{%
  \reset@font\fontsize{#1}{#2pt}%
  \fontfamily{#3}\fontseries{#4}\fontshape{#5}%
  \selectfont}%
\fi\endgroup%
\begin{picture}(1872,5942)(841,-8016)
\put(856,-4201){\makebox(0,0)[rb]{\smash{{\SetFigFontNFSS{20}{24.0}{\rmdefault}{\mddefault}{\updefault}{\color[rgb]{0,0,0}$\delta\mathord+\omega$}%
}}}}
\put(856,-3526){\makebox(0,0)[rb]{\smash{{\SetFigFontNFSS{20}{24.0}{\rmdefault}{\mddefault}{\updefault}{\color[rgb]{0,0,0}$\delta\mathord+\omega\mathord\cdot 2$}%
}}}}
\put(856,-4876){\makebox(0,0)[rb]{\smash{{\SetFigFontNFSS{20}{24.0}{\rmdefault}{\mddefault}{\updefault}{\color[rgb]{0,0,0}$\delta$}%
}}}}
\put(1576,-7521){\makebox(0,0)[b]{\smash{{\SetFigFontNFSS{20}{24.0}{\rmdefault}{\mddefault}{\updefault}{\color[rgb]{0,0,0}$0$}%
}}}}
\put(1136,-7521){\makebox(0,0)[b]{\smash{{\SetFigFontNFSS{20}{24.0}{\rmdefault}{\mddefault}{\updefault}{\color[rgb]{0,0,0}$\sim$}%
}}}}
\put(2046,-7521){\makebox(0,0)[b]{\smash{{\SetFigFontNFSS{20}{24.0}{\rmdefault}{\mddefault}{\updefault}{\color[rgb]{0,0,0}$1$}%
}}}}
\put(2476,-7521){\makebox(0,0)[b]{\smash{{\SetFigFontNFSS{20}{24.0}{\rmdefault}{\mddefault}{\updefault}{\color[rgb]{0,0,0}$2$}%
}}}}
\put(1801,-7911){\makebox(0,0)[b]{\smash{{\SetFigFontNFSS{20}{24.0}{\rmdefault}{\mddefault}{\updefault}{\color[rgb]{0,0,0}\mylabeltypicalfinite}%
}}}}
\end{picture}%

%% file: backandforth.pstex_t
\begin{picture}(0,0)%
\includegraphics{backandforth.pstex}%
\end{picture}%
\setlength{\unitlength}{4144sp}%
\begingroup\makeatletter\ifx\SetFigFontNFSS\undefined%
\gdef\SetFigFontNFSS#1#2#3#4#5{%
  \reset@font\fontsize{#1}{#2pt}%
  \fontfamily{#3}\fontseries{#4}\fontshape{#5}%
  \selectfont}%
\fi\endgroup%
\begin{picture}(2547,6167)(841,-8061)
\put(991,-7486){\makebox(0,0)[b]{\smash{{\SetFigFontNFSS{20}{24.0}{\rmdefault}{\mddefault}{\updefault}{\color[rgb]{0,0,0}$\sim$}%
}}}}
\put(2156,-7911){\makebox(0,0)[b]{\smash{{\SetFigFontNFSS{20}{24.0}{\rmdefault}{\mddefault}{\updefault}{\color[rgb]{0,0,0}\mylabelbackandforth}%
}}}}
\put(2144,-7486){\makebox(0,0)[b]{\smash{{\SetFigFontNFSS{20}{24.0}{\rmdefault}{\mddefault}{\updefault}{\color[rgb]{0,0,0}$n$}%
}}}}
\put(856,-4876){\makebox(0,0)[rb]{\smash{{\SetFigFontNFSS{20}{24.0}{\rmdefault}{\mddefault}{\updefault}{\color[rgb]{0,0,0}$\beta\mathord+\omega$}%
}}}}
\put(856,-5776){\makebox(0,0)[rb]{\smash{{\SetFigFontNFSS{20}{24.0}{\rmdefault}{\mddefault}{\updefault}{\color[rgb]{0,0,0}$\beta$}%
}}}}
\put(856,-3076){\makebox(0,0)[rb]{\smash{{\SetFigFontNFSS{20}{24.0}{\rmdefault}{\mddefault}{\updefault}{\color[rgb]{0,0,0}$\delta$}%
}}}}
\put(856,-2176){\makebox(0,0)[rb]{\smash{{\SetFigFontNFSS{20}{24.0}{\rmdefault}{\mddefault}{\updefault}{\color[rgb]{0,0,0}$\delta\mathord+\omega$}%
}}}}
\end{picture}%

%% file: stacking.pstex_t
\begin{picture}(0,0)%
\includegraphics{stacking.pstex}%
\end{picture}%
\setlength{\unitlength}{4144sp}%
\begingroup\makeatletter\ifx\SetFigFontNFSS\undefined%
\gdef\SetFigFontNFSS#1#2#3#4#5{%
  \reset@font\fontsize{#1}{#2pt}%
  \fontfamily{#3}\fontseries{#4}\fontshape{#5}%
  \selectfont}%
\fi\endgroup%
\begin{picture}(5247,6357)(-1859,-7996)
\put(181,-6451){\makebox(0,0)[rb]{\smash{{\SetFigFontNFSS{20}{24.0}{\rmdefault}{\mddefault}{\updefault}{\color[rgb]{0,0,0}$\omega$}%
}}}}
\put(-114,-7891){\makebox(0,0)[b]{\smash{{\SetFigFontNFSS{20}{24.0}{\rmdefault}{\mddefault}{\updefault}{\color[rgb]{0,0,0}\mylabelstacking}%
}}}}
\put(181,-4201){\makebox(0,0)[rb]{\smash{{\SetFigFontNFSS{20}{24.0}{\rmdefault}{\mddefault}{\updefault}{\color[rgb]{0,0,0}$\delta\mathord+\omega$}%
}}}}
\put(181,-5101){\makebox(0,0)[rb]{\smash{{\SetFigFontNFSS{20}{24.0}{\rmdefault}{\mddefault}{\updefault}{\color[rgb]{0,0,0}$\delta$}%
}}}}
\put(2386,-7531){\makebox(0,0)[b]{\smash{{\SetFigFontNFSS{20}{24.0}{\rmdefault}{\mddefault}{\updefault}{\color[rgb]{0,0,0}$\sim$}%
}}}}
\put(361,-7531){\makebox(0,0)[b]{\smash{{\SetFigFontNFSS{20}{24.0}{\rmdefault}{\mddefault}{\updefault}{\color[rgb]{0,0,0}$\sim$}%
}}}}
\put(2206,-6451){\makebox(0,0)[rb]{\smash{{\SetFigFontNFSS{20}{24.0}{\rmdefault}{\mddefault}{\updefault}{\color[rgb]{0,0,0}$\omega$}%
}}}}
\put(2206,-5101){\makebox(0,0)[rb]{\smash{{\SetFigFontNFSS{20}{24.0}{\rmdefault}{\mddefault}{\updefault}{\color[rgb]{0,0,0}$\delta$}%
}}}}
\put(2206,-4201){\makebox(0,0)[rb]{\smash{{\SetFigFontNFSS{20}{24.0}{\rmdefault}{\mddefault}{\updefault}{\color[rgb]{0,0,0}$\delta\mathord+\omega$}%
}}}}
\put(-1844,-5101){\makebox(0,0)[rb]{\smash{{\SetFigFontNFSS{20}{24.0}{\rmdefault}{\mddefault}{\updefault}{\color[rgb]{0,0,0}$0$}%
}}}}
\put(-1844,-2176){\makebox(0,0)[rb]{\smash{{\SetFigFontNFSS{20}{24.0}{\rmdefault}{\mddefault}{\updefault}{\color[rgb]{0,0,0}$\delta'\mathord+\omega$}%
}}}}
\put(-1844,-4201){\makebox(0,0)[rb]{\smash{{\SetFigFontNFSS{20}{24.0}{\rmdefault}{\mddefault}{\updefault}{\color[rgb]{0,0,0}$\omega$}%
}}}}
\put(-1844,-3076){\makebox(0,0)[rb]{\smash{{\SetFigFontNFSS{20}{24.0}{\rmdefault}{\mddefault}{\updefault}{\color[rgb]{0,0,0}$\delta'$}%
}}}}
\put(2206,-3076){\makebox(0,0)[rb]{\smash{{\SetFigFontNFSS{20}{24.0}{\rmdefault}{\mddefault}{\updefault}{\color[rgb]{0,0,0}$\delta\mathord+\delta'$}%
}}}}
\put(2206,-2176){\makebox(0,0)[rb]{\smash{{\SetFigFontNFSS{20}{24.0}{\rmdefault}{\mddefault}{\updefault}{\color[rgb]{0,0,0}$\delta\mathord+\delta'\mathord+\omega$}%
}}}}
\put(2611,-2626){\makebox(0,0)[b]{\smash{{\SetFigFontNFSS{20}{24.0}{\rmdefault}{\mddefault}{\updefault}{\color[rgb]{0,0,0}$\bar y$}%
}}}}
\put(-1439,-2626){\makebox(0,0)[b]{\smash{{\SetFigFontNFSS{20}{24.0}{\rmdefault}{\mddefault}{\updefault}{\color[rgb]{0,0,0}$\bar y$}%
}}}}
\put(721,-4876){\makebox(0,0)[b]{\smash{{\SetFigFontNFSS{20}{24.0}{\rmdefault}{\mddefault}{\updefault}{\color[rgb]{0,0,0}$\bar x$}%
}}}}
\put(-1259,-1906){\makebox(0,0)[b]{\smash{{\SetFigFontNFSS{20}{24.0}{\rmdefault}{\mddefault}{\updefault}{\color[rgb]{0,0,0}$q'$}%
}}}}
\put(811,-3931){\makebox(0,0)[b]{\smash{{\SetFigFontNFSS{20}{24.0}{\rmdefault}{\mddefault}{\updefault}{\color[rgb]{0,0,0}$p$}%
}}}}
\put(2791,-1951){\makebox(0,0)[b]{\smash{{\SetFigFontNFSS{20}{24.0}{\rmdefault}{\mddefault}{\updefault}{\color[rgb]{0,0,0}$q=p\Lsh q'$}%
}}}}
\put(-1124,-5247){\makebox(0,0)[b]{\smash{{\SetFigFontNFSS{20}{24.0}{\rmdefault}{\mddefault}{\updefault}{\color[rgb]{0,0,0}$n$}%
}}}}
\put(901,-7531){\makebox(0,0)[b]{\smash{{\SetFigFontNFSS{20}{24.0}{\rmdefault}{\mddefault}{\updefault}{\color[rgb]{0,0,0}$n$}%
}}}}
\put(2926,-7531){\makebox(0,0)[b]{\smash{{\SetFigFontNFSS{20}{24.0}{\rmdefault}{\mddefault}{\updefault}{\color[rgb]{0,0,0}$n$}%
}}}}
\put(1756,-5776){\makebox(0,0)[b]{\smash{{\SetFigFontNFSS{20}{24.0}{\rmdefault}{\mddefault}{\updefault}{\color[rgb]{0,0,0}$\phi$}%
}}}}
\put(811,-3436){\makebox(0,0)[b]{\smash{{\SetFigFontNFSS{20}{24.0}{\rmdefault}{\mddefault}{\updefault}{\color[rgb]{0,0,0}$\text{Id}$}%
}}}}
\put(-449,-4651){\makebox(0,0)[lb]{\smash{{\SetFigFontNFSS{20}{24.0}{\rmdefault}{\mddefault}{\updefault}{\color[rgb]{0,0,0}$\phi$}%
}}}}
\end{picture}%

%% file: strangea.pstex_t
\begin{picture}(0,0)%
\includegraphics{strangea.pstex}%
\end{picture}%
\setlength{\unitlength}{4144sp}%
\begingroup\makeatletter\ifx\SetFigFontNFSS\undefined%
\gdef\SetFigFontNFSS#1#2#3#4#5{%
  \reset@font\fontsize{#1}{#2pt}%
  \fontfamily{#3}\fontseries{#4}\fontshape{#5}%
  \selectfont}%
\fi\endgroup%
\begin{picture}(4299,6372)(-11,-11461)
\put(3712,-10951){\makebox(0,0)[b]{\smash{{\SetFigFontNFSS{20}{24.0}{\rmdefault}{\mddefault}{\updefault}{\color[rgb]{0,0,0}$2$}%
}}}}
\put(3488,-10951){\makebox(0,0)[b]{\smash{{\SetFigFontNFSS{20}{24.0}{\rmdefault}{\mddefault}{\updefault}{\color[rgb]{0,0,0}$1$}%
}}}}
\put(3263,-10951){\makebox(0,0)[b]{\smash{{\SetFigFontNFSS{20}{24.0}{\rmdefault}{\mddefault}{\updefault}{\color[rgb]{0,0,0}$0$}%
}}}}
\put(1688,-10951){\makebox(0,0)[b]{\smash{{\SetFigFontNFSS{20}{24.0}{\rmdefault}{\mddefault}{\updefault}{\color[rgb]{0,0,0}$0$}%
}}}}
\put(1913,-10951){\makebox(0,0)[b]{\smash{{\SetFigFontNFSS{20}{24.0}{\rmdefault}{\mddefault}{\updefault}{\color[rgb]{0,0,0}$1$}%
}}}}
\put(2137,-10951){\makebox(0,0)[b]{\smash{{\SetFigFontNFSS{20}{24.0}{\rmdefault}{\mddefault}{\updefault}{\color[rgb]{0,0,0}$2$}%
}}}}
\put(562,-10951){\makebox(0,0)[b]{\smash{{\SetFigFontNFSS{20}{24.0}{\rmdefault}{\mddefault}{\updefault}{\color[rgb]{0,0,0}$2$}%
}}}}
\put(338,-10951){\makebox(0,0)[b]{\smash{{\SetFigFontNFSS{20}{24.0}{\rmdefault}{\mddefault}{\updefault}{\color[rgb]{0,0,0}$1$}%
}}}}
\put(113,-10951){\makebox(0,0)[b]{\smash{{\SetFigFontNFSS{20}{24.0}{\rmdefault}{\mddefault}{\updefault}{\color[rgb]{0,0,0}$0$}%
}}}}
\put(2161,-11311){\makebox(0,0)[b]{\smash{{\SetFigFontNFSS{20}{24.0}{\rmdefault}{\mddefault}{\updefault}{\color[rgb]{0,0,0}\labelstrangea}%
}}}}
\put(3691,-5371){\makebox(0,0)[b]{\smash{{\SetFigFontNFSS{20}{24.0}{\rmdefault}{\mddefault}{\updefault}{\color[rgb]{0,0,0}$q$}%
}}}}
\put(2116,-6721){\makebox(0,0)[b]{\smash{{\SetFigFontNFSS{20}{24.0}{\rmdefault}{\mddefault}{\updefault}{\color[rgb]{0,0,0}$p_2$}%
}}}}
\put(541,-8251){\makebox(0,0)[b]{\smash{{\SetFigFontNFSS{20}{24.0}{\rmdefault}{\mddefault}{\updefault}{\color[rgb]{0,0,0}$p_1$}%
}}}}
\end{picture}%

%% file: strangeb.pstex_t
\begin{picture}(0,0)%
\includegraphics{strangeb.pstex}%
\end{picture}%
\setlength{\unitlength}{4144sp}%
\begingroup\makeatletter\ifx\SetFigFontNFSS\undefined%
\gdef\SetFigFontNFSS#1#2#3#4#5{%
  \reset@font\fontsize{#1}{#2pt}%
  \fontfamily{#3}\fontseries{#4}\fontshape{#5}%
  \selectfont}%
\fi\endgroup%
\begin{picture}(4299,6372)(-11,-11461)
\put(3714,-10950){\makebox(0,0)[b]{\smash{{\SetFigFontNFSS{20}{24.0}{\rmdefault}{\mddefault}{\updefault}{\color[rgb]{0,0,0}$2$}%
}}}}
\put(3490,-10950){\makebox(0,0)[b]{\smash{{\SetFigFontNFSS{20}{24.0}{\rmdefault}{\mddefault}{\updefault}{\color[rgb]{0,0,0}$1$}%
}}}}
\put(3265,-10950){\makebox(0,0)[b]{\smash{{\SetFigFontNFSS{20}{24.0}{\rmdefault}{\mddefault}{\updefault}{\color[rgb]{0,0,0}$0$}%
}}}}
\put(1690,-10950){\makebox(0,0)[b]{\smash{{\SetFigFontNFSS{20}{24.0}{\rmdefault}{\mddefault}{\updefault}{\color[rgb]{0,0,0}$0$}%
}}}}
\put(1915,-10950){\makebox(0,0)[b]{\smash{{\SetFigFontNFSS{20}{24.0}{\rmdefault}{\mddefault}{\updefault}{\color[rgb]{0,0,0}$1$}%
}}}}
\put(2139,-10950){\makebox(0,0)[b]{\smash{{\SetFigFontNFSS{20}{24.0}{\rmdefault}{\mddefault}{\updefault}{\color[rgb]{0,0,0}$2$}%
}}}}
\put(564,-10950){\makebox(0,0)[b]{\smash{{\SetFigFontNFSS{20}{24.0}{\rmdefault}{\mddefault}{\updefault}{\color[rgb]{0,0,0}$2$}%
}}}}
\put(340,-10950){\makebox(0,0)[b]{\smash{{\SetFigFontNFSS{20}{24.0}{\rmdefault}{\mddefault}{\updefault}{\color[rgb]{0,0,0}$1$}%
}}}}
\put(115,-10950){\makebox(0,0)[b]{\smash{{\SetFigFontNFSS{20}{24.0}{\rmdefault}{\mddefault}{\updefault}{\color[rgb]{0,0,0}$0$}%
}}}}
\put(541,-10591){\makebox(0,0)[b]{\smash{{\SetFigFontNFSS{20}{24.0}{\rmdefault}{\mddefault}{\updefault}{\color[rgb]{0,0,0}$0$}%
}}}}
\put(2161,-10591){\makebox(0,0)[b]{\smash{{\SetFigFontNFSS{20}{24.0}{\rmdefault}{\mddefault}{\updefault}{\color[rgb]{0,0,0}$x$}%
}}}}
\put(1891,-9241){\makebox(0,0)[b]{\smash{{\SetFigFontNFSS{20}{24.0}{\rmdefault}{\mddefault}{\updefault}{\color[rgb]{0,0,0}$x$}%
}}}}
\put(1666,-7666){\makebox(0,0)[b]{\smash{{\SetFigFontNFSS{20}{24.0}{\rmdefault}{\mddefault}{\updefault}{\color[rgb]{0,0,0}$x$}%
}}}}
\put(3736,-10591){\makebox(0,0)[b]{\smash{{\SetFigFontNFSS{20}{24.0}{\rmdefault}{\mddefault}{\updefault}{\color[rgb]{0,0,0}$0$}%
}}}}
\put(3511,-9241){\makebox(0,0)[b]{\smash{{\SetFigFontNFSS{20}{24.0}{\rmdefault}{\mddefault}{\updefault}{\color[rgb]{0,0,0}$0$}%
}}}}
\put(3241,-7666){\makebox(0,0)[b]{\smash{{\SetFigFontNFSS{20}{24.0}{\rmdefault}{\mddefault}{\updefault}{\color[rgb]{0,0,0}$0$}%
}}}}
\put(2161,-11311){\makebox(0,0)[b]{\smash{{\SetFigFontNFSS{20}{24.0}{\rmdefault}{\mddefault}{\updefault}{\color[rgb]{0,0,0}\labelstrangeb}%
}}}}
\put(3691,-5371){\makebox(0,0)[b]{\smash{{\SetFigFontNFSS{20}{24.0}{\rmdefault}{\mddefault}{\updefault}{\color[rgb]{0,0,0}$q$}%
}}}}
\put(2116,-6721){\makebox(0,0)[b]{\smash{{\SetFigFontNFSS{20}{24.0}{\rmdefault}{\mddefault}{\updefault}{\color[rgb]{0,0,0}$p_2$}%
}}}}
\put(541,-8251){\makebox(0,0)[b]{\smash{{\SetFigFontNFSS{20}{24.0}{\rmdefault}{\mddefault}{\updefault}{\color[rgb]{0,0,0}$p_1$}%
}}}}
\end{picture}%